\parindent=12pt   \parskip=0pt     
\pageno=1 


\hoffset=15mm    
\voffset=1cm    
 

\ifnum\mag=\magstep1
\hoffset=-2mm   
\voffset=.8cm   
\fi


\pretolerance=500 \tolerance=1000  \brokenpenalty=5000

\catcode`\@=11

\font\eightrm=cmr8         \font\eighti=cmmi8
\font\eightsy=cmsy8        \font\eightbf=cmbx8
\font\eighttt=cmtt8        \font\eightit=cmti8
\font\eightsl=cmsl8        \font\sixrm=cmr6
\font\sixi=cmmi6           \font\sixsy=cmsy6
\font\sixbf=cmbx6


\font\tengoth=eufm10       \font\tenbboard=msbm10
\font\eightgoth=eufm10 at 8pt      \font\eightbboard=msbm10 at 8 pt
\font\sevengoth=eufm7      \font\sevenbboard=msbm7
\font\sixgoth=eufm7 at 6 pt        \font\fivegoth=eufm5

 \font\tencyr=wncyr10       
\font\eightcyr=wncyr10 at 8 pt      
\font\sevencyr=wncyr10 at 7 pt      
\font\sixcyr=wncyr10 at 6 pt


\skewchar\eighti='177 \skewchar\sixi='177
\skewchar\eightsy='60 \skewchar\sixsy='60


\newfam\gothfam           \newfam\bboardfam
\newfam\cyrfam

\def\tenpoint{%
  \textfont0=\tenrm \scriptfont0=\sevenrm \scriptscriptfont0=\fiverm
  \def\rm{\fam\z@\tenrm}%
  \textfont1=\teni  \scriptfont1=\seveni  \scriptscriptfont1=\fivei
  \def\oldstyle{\fam\@ne\teni}\let\old=\oldstyle
  \textfont2=\tensy \scriptfont2=\sevensy \scriptscriptfont2=\fivesy
  \textfont\gothfam=\tengoth \scriptfont\gothfam=\sevengoth
  \scriptscriptfont\gothfam=\fivegoth
  \def\goth{\fam\gothfam\tengoth}%
  \textfont\bboardfam=\tenbboard \scriptfont\bboardfam=\sevenbboard
  \scriptscriptfont\bboardfam=\sevenbboard
  \def\bb{\fam\bboardfam\tenbboard}%
 \textfont\cyrfam=\tencyr \scriptfont\cyrfam=\sevencyr
  \scriptscriptfont\cyrfam=\sixcyr
  \def\cyr{\fam\cyrfam\tencyr}%
  \textfont\itfam=\tenit
  \def\it{\fam\itfam\tenit}%
  \textfont\slfam=\tensl
  \def\sl{\fam\slfam\tensl}%
  \textfont\bffam=\tenbf \scriptfont\bffam=\sevenbf
  \scriptscriptfont\bffam=\fivebf
  \def\bf{\fam\bffam\tenbf}%
  \textfont\ttfam=\tentt
  \def\tt{\fam\ttfam\tentt}%
  \abovedisplayskip=12pt plus 3pt minus 9pt
  \belowdisplayskip=\abovedisplayskip
  \abovedisplayshortskip=0pt plus 3pt
  \belowdisplayshortskip=4pt plus 3pt 
  \smallskipamount=3pt plus 1pt minus 1pt
  \medskipamount=6pt plus 2pt minus 2pt
  \bigskipamount=12pt plus 4pt minus 4pt
  \normalbaselineskip=12pt
  \setbox\strutbox=\hbox{\vrule height8.5pt depth3.5pt width0pt}%
  \let\bigf@nt=\tenrm       \let\smallf@nt=\sevenrm
  \normalbaselines\rm}

\def\eightpoint{%
  \textfont0=\eightrm \scriptfont0=\sixrm \scriptscriptfont0=\fiverm
  \def\rm{\fam\z@\eightrm}%
  \textfont1=\eighti  \scriptfont1=\sixi  \scriptscriptfont1=\fivei
  \def\oldstyle{\fam\@ne\eighti}\let\old=\oldstyle
  \textfont2=\eightsy \scriptfont2=\sixsy \scriptscriptfont2=\fivesy
  \textfont\gothfam=\eightgoth \scriptfont\gothfam=\sixgoth
  \scriptscriptfont\gothfam=\fivegoth
  \def\goth{\fam\gothfam\eightgoth}%
  \textfont\cyrfam=\eightcyr \scriptfont\cyrfam=\sixcyr
  \scriptscriptfont\cyrfam=\sixcyr
  \def\cyr{\fam\cyrfam\eightcyr}%
  \textfont\bboardfam=\eightbboard \scriptfont\bboardfam=\sevenbboard
  \scriptscriptfont\bboardfam=\sevenbboard
  \def\bb{\fam\bboardfam}%
  \textfont\itfam=\eightit
  \def\it{\fam\itfam\eightit}%
  \textfont\slfam=\eightsl
  \def\sl{\fam\slfam\eightsl}%
  \textfont\bffam=\eightbf \scriptfont\bffam=\sixbf
  \scriptscriptfont\bffam=\fivebf
  \def\bf{\fam\bffam\eightbf}%
  \textfont\ttfam=\eighttt
  \def\tt{\fam\ttfam\eighttt}%
  \abovedisplayskip=9pt plus 3pt minus 9pt
  \belowdisplayskip=\abovedisplayskip
  \abovedisplayshortskip=0pt plus 3pt
  \belowdisplayshortskip=3pt plus 3pt 
  \smallskipamount=2pt plus 1pt minus 1pt
  \medskipamount=4pt plus 2pt minus 1pt
  \bigskipamount=9pt plus 3pt minus 3pt
  \normalbaselineskip=9pt
  \setbox\strutbox=\hbox{\vrule height7pt depth2pt width0pt}%
  \let\bigf@nt=\eightrm     \let\smallf@nt=\sixrm
  \normalbaselines\rm}

\tenpoint


\def\pc#1{\bigf@nt#1\smallf@nt}         \def\pd#1 {{\pc#1} }


\catcode`\;=\active
\def;{\relax\ifhmode\ifdim\lastskip>\z@\unskip\fi
\kern\fontdimen2  -1.2 \fontdimen3 \string;}

\catcode`\:=\active
\def:{\relax\ifhmode\ifdim\lastskip>\z@\unskip\fi\penalty\@M\ \fi\string:}

\catcode`\!=\active
\def!{\relax\ifhmode\ifdim\lastskip>\z@
\unskip\fi\kern\fontdimen2  -1.1 \fontdimen3 \string!}

\catcode`\?=\active
\def?{\relax\ifhmode\ifdim\lastskip>\z@
\unskip\fi\kern\fontdimen2  -1.1 \fontdimen3 \string?}

\def\^#1{\if#1i{\accent"5E\i}\else{\accent"5E #1}\fi}
\def\"#1{\if#1i{\accent"7F\i}\else{\accent"7F #1}\fi}

\frenchspacing


\newtoks\auteurcourant      \auteurcourant={\hfil}
\newtoks\titrecourant       \titrecourant={\hfil}

\newtoks\hautpagetitre      \hautpagetitre={\hfil}
\newtoks\baspagetitre       \baspagetitre={\hfil}

\newtoks\hautpagegauche     
\hautpagegauche={\eightpoint\rlap{\folio}\hfil\the\auteurcourant\hfil}
\newtoks\hautpagedroite     
\hautpagedroite={\eightpoint\hfil\the\titrecourant\hfil\llap{\folio}}

\newtoks\baspagegauche      \baspagegauche={\hfil} 
\newtoks\baspagedroite      \baspagedroite={\hfil}

\newif\ifpagetitre          \pagetitretrue  


\headline={\ifpagetitre\the\hautpagetitre
\else\ifodd\pageno\the\hautpagedroite\else\the\hautpagegauche\fi\fi}

\footline={\ifpagetitre\the\baspagetitre\else
\ifodd\pageno\the\baspagedroite\else\the\baspagegauche\fi\fi
\global\pagetitrefalse}


\def\raggedbottom{\topskip 10pt plus 36pt\r@ggedbottomtrue}



\def\pointir{\unskip . --- \ignorespaces}


\def\Bigbreak{\vskip-\lastskip\bigbreak}
\def\Medbreak{\vskip-\lastskip\medbreak}


\def\ctexte#1\endctexte{%
  \hbox{$\vcenter{\halign{\hfill##\hfill\crcr#1\crcr}}$}}


\long\def\ctitre#1\endctitre{%
    \ifdim\lastskip<24pt\vskip-\lastskip\bigbreak\bigbreak\fi
  		\vbox{\parindent=0pt\leftskip=0pt plus 1fill
          \rightskip=\leftskip
          \parfillskip=0pt\bf#1\par}
    \bigskip\nobreak}

\long\def\section#1\endsection{%
\vskip 0pt plus 3\normalbaselineskip
\penalty-250
\vskip 0pt plus -3\normalbaselineskip
\Bigbreak
\message{[section \string: #1]}{\bf#1\unskip}\pointir}

\long\def\sectiona#1\endsection{%
\vskip 0pt plus 3\normalbaselineskip
\penalty-250
\vskip 0pt plus -3\normalbaselineskip
\Bigbreak
\message{[sectiona \string: #1]}%
{\bf#1}\medskip\nobreak}

\long\def\subsection#1\endsubsection{%
\Medbreak
{\it#1\unskip}\pointir}

\long\def\subsectiona#1\endsubsection{%
\Medbreak
{\it#1}\par\nobreak}

\def\rem#1\endrem{%
\Medbreak
{\it#1\unskip} : }

\def\remp#1\endrem{%
\Medbreak
{\pc #1\unskip}\pointir}

\def\rema#1\endrem{%
\Medbreak
{\it #1}\par\nobreak}

\def\newparwithcolon#1\endnewparwithcolon{
\Medbreak
{#1\unskip} : }

\def\newparwithpointir#1\endnewparwithpointir{
\Medbreak
{#1\unskip}\pointir}

\def\newpara#1\endnewpar{
\Medbreak
{#1\unskip}\smallskip\nobreak}


\long\def\th#1 #2\enonce#3\endth{%
   \Medbreak
   {\pc#1} {#2\unskip}\pointir{\it #3}\medskip}

\long\def\tha#1 #2\enonce#3\endth{%
   \Medbreak
   {\pc#1} {#2\unskip}\par\nobreak{\it #3}\medskip}


\long\def\Th#1 #2 #3\enonce#4\endth{%
   \Medbreak
   #1 {\pc#2} {#3\unskip}\pointir{\it #4}\medskip}

\long\def\Tha#1 #2 #3\enonce#4\endth{%
   \Medbreak
   #1 {\pc#2} #3\par\nobreak{\it #4}\medskip}


\def\decale#1{\smallbreak\hskip 28pt\llap{#1}\kern 5pt}
\def\decaledecale#1{\smallbreak\hskip 34pt\llap{#1}\kern 5pt}
\def\puce{\smallbreak\hskip 6pt{$\scriptstyle\bullet$}\kern 5pt}



\def\displaylinesno#1{\displ@y\halign{
\hbox to\displaywidth{$\@lign\hfil\displaystyle##\hfil$}&
\llap{$##$}\crcr#1\crcr}}


\def\ldisplaylinesno#1{\displ@y\halign{ 
\hbox to\displaywidth{$\@lign\hfil\displaystyle##\hfil$}&
\kern-\displaywidth\rlap{$##$}\tabskip\displaywidth\crcr#1\crcr}}


\def\eqalign#1{\null\,\vcenter{\openup\jot\m@th\ialign{
\strut\hfil$\displaystyle{##}$&$\displaystyle{{}##}$\hfil
&&\quad\strut\hfil$\displaystyle{##}$&$\displaystyle{{}##}$\hfil
\crcr#1\crcr}}\,}


\def\system#1{\left\{\null\,\vcenter{\openup1\jot\m@th
\ialign{\strut$##$&\hfil$##$&$##$\hfil&&
        \enskip$##$\enskip&\hfil$##$&$##$\hfil\crcr#1\crcr}}\right.}


\let\@ldmessage=\message

\def\message#1{{\def\pc{\string\pc\space}%
                \def\'{\string'}\def\`{\string`}%
                \def\^{\string^}\def\"{\string"}%
                \@ldmessage{#1}}}



\def\up#1{\raise 1ex\hbox{\smallf@nt#1}}


\def\qed{\raise -2pt\hbox{\vrule\vbox to 10pt{\hrule width 4pt
                 \vfill\hrule}\vrule}}

\def\virg{\raise .4ex\hbox{,}}   


\def\build#1_#2^#3{\mathrel{
\mathop{\kern 0pt#1}\limits_{#2}^{#3}}}


\def\boxit#1#2{%
\setbox1=\hbox{\kern#1{#2}\kern#1}%
\dimen1=\ht1 \advance\dimen1 by #1 \dimen2=\dp1 \advance\dimen2 by #1 
\setbox1=\hbox{\vrule height\dimen1 depth\dimen2\box1\vrule}%
\setbox1=\vbox{\hrule\box1\hrule}%
\advance\dimen1 by .6pt \ht1=\dimen1 
\advance\dimen2 by .6pt \dp1=\dimen2  \box1\relax}


\catcode`\@=12

\showboxbreadth=-1  \showboxdepth=-1




\def\Grille{\setbox13=\vbox to 5\unitlength{\hrule width 109mm\vfill} 
\setbox13=\vbox to 65mm{\offinterlineskip\leaders\copy13\vfill\kern-1pt\hrule} 
\setbox14=\hbox to 5\unitlength{\vrule height 65mm\hfill} 
\setbox14=\hbox to 109mm{\leaders\copy14\hfill\kern-2mm\vrule height 65mm}
\ht14=0pt\dp14=0pt\wd14=0pt \setbox13=\vbox to
0pt{\vss\box13\offinterlineskip\box14} \wd13=0pt\box13}


\def\fleche(#1,#2)\dir(#3,#4)\long#5{%
\noalign{\leftput(#1,#2){\vector(#3,#4){#5}}}}

\def\ligne(#1,#2)\dir(#3,#4)\long#5{%
\noalign{\leftput(#1,#2){\lline(#3,#4){#5}}}}

\def\put(#1,#2)#3{\noalign{\setbox1=\hbox{%
    \kern #1\unitlength
    \raise #2\unitlength\hbox{$#3$}}%
    \ht1=0pt \wd1=0pt \dp1=0pt\box1}}


\def\diagram#1{\def\normalbaselines{\baselineskip=0pt\lineskip=5pt}
\matrix{#1}}

\def\hfl#1#2#3{\smash{\mathop{\hbox to#3{\rightarrowfill}}\limits
^{\scriptstyle#1}_{\scriptstyle#2}}}

\def\gfl#1#2#3{\smash{\mathop{\hbox to#3{\leftarrowfill}}\limits
^{\scriptstyle#1}_{\scriptstyle#2}}}

\def\vfl#1#2#3{\llap{$\scriptstyle #1$}
\left\downarrow\vbox to#3{}\right.\rlap{$\scriptstyle #2$}}


 \message{`lline' & `vector' macros from LaTeX}
 \catcode`@=11
\def\{{\relax\ifmmode\lbrace\else$\lbrace$\fi}
\def\}{\relax\ifmmode\rbrace\else$\rbrace$\fi}
\def\newcount{\alloc@0\count\countdef\insc@unt}
\def\newdimen{\alloc@1\dimen\dimendef\insc@unt}
\def\newwrite{\alloc@7\write\chardef\sixt@@n}

\newwrite\@unused
\newcount\@tempcnta
\newcount\@tempcntb
\newdimen\@tempdima
\newdimen\@tempdimb
\newbox\@tempboxa

\def\@spaces{\space\space\space\space}
\def\@whilenoop#1{}
\def\@whiledim#1\do #2{\ifdim #1\relax#2\@iwhiledim{#1\relax#2}\fi}
\def\@iwhiledim#1{\ifdim #1\let\@nextwhile=\@iwhiledim
        \else\let\@nextwhile=\@whilenoop\fi\@nextwhile{#1}}
\def\@badlinearg{\@latexerr{Bad \string\line\space or \string\vector
   \space argument}}
\def\@latexerr#1#2{\begingroup
\edef\@tempc{#2}\expandafter\errhelp\expandafter{\@tempc}%
\def\@eha{Your command was ignored.
^^JType \space I <command> <return> \space to replace it
  with another command,^^Jor \space <return> \space to continue without
it.} 
\def\@ehb{You've lost some text. \space \@ehc}
\def\@ehc{Try typing \space <return>
  \space to proceed.^^JIf that doesn't work, type \space X <return> \space to
  quit.}
\def\@ehd{You're in trouble here.  \space\@ehc}

\typeout{LaTeX error. \space See LaTeX manual for explanation.^^J
 \space\@spaces\@spaces\@spaces Type \space H <return> \space for
 immediate help.}\errmessage{#1}\endgroup}
\def\typeout#1{{\let\protect\string\immediate\write\@unused{#1}}}

\font\tenln    = line10
\font\tenlnw   = linew10

\newdimen\@wholewidth
\newdimen\@halfwidth
\newdimen\unitlength 

\unitlength =1pt


\def\thinlines{\let\@linefnt\tenln \let\@circlefnt\tencirc
  \@wholewidth\fontdimen8\tenln \@halfwidth .5\@wholewidth}
\def\thicklines{\let\@linefnt\tenlnw \let\@circlefnt\tencircw
  \@wholewidth\fontdimen8\tenlnw \@halfwidth .5\@wholewidth}

\def\linethickness#1{\@wholewidth #1\relax \@halfwidth .5\@wholewidth}

\newif\if@negarg

\def\lline(#1,#2)#3{\@xarg #1\relax \@yarg #2\relax
\@linelen=#3\unitlength
\ifnum\@xarg =0 \@vline
  \else \ifnum\@yarg =0 \@hline \else \@sline\fi
\fi}

\def\@sline{\ifnum\@xarg< 0 \@negargtrue \@xarg -\@xarg \@yyarg -\@yarg
  \else \@negargfalse \@yyarg \@yarg \fi
\ifnum \@yyarg >0 \@tempcnta\@yyarg \else \@tempcnta -\@yyarg \fi
\ifnum\@tempcnta>6 \@badlinearg\@tempcnta0 \fi
\setbox\@linechar\hbox{\@linefnt\@getlinechar(\@xarg,\@yyarg)}%
\ifnum \@yarg >0 \let\@upordown\raise \@clnht\z@
   \else\let\@upordown\lower \@clnht \ht\@linechar\fi
\@clnwd=\wd\@linechar
\if@negarg \hskip -\wd\@linechar \def\@tempa{\hskip -2\wd\@linechar}\else
     \let\@tempa\relax \fi
\@whiledim \@clnwd <\@linelen \do
  {\@upordown\@clnht\copy\@linechar
   \@tempa
   \advance\@clnht \ht\@linechar
   \advance\@clnwd \wd\@linechar}%
\advance\@clnht -\ht\@linechar
\advance\@clnwd -\wd\@linechar
\@tempdima\@linelen\advance\@tempdima -\@clnwd
\@tempdimb\@tempdima\advance\@tempdimb -\wd\@linechar
\if@negarg \hskip -\@tempdimb \else \hskip \@tempdimb \fi
\multiply\@tempdima \@m
\@tempcnta \@tempdima \@tempdima \wd\@linechar \divide\@tempcnta \@tempdima
\@tempdima \ht\@linechar \multiply\@tempdima \@tempcnta
\divide\@tempdima \@m
\advance\@clnht \@tempdima
\ifdim \@linelen <\wd\@linechar
   \hskip \wd\@linechar
  \else\@upordown\@clnht\copy\@linechar\fi}

\def\@hline{\ifnum \@xarg <0 \hskip -\@linelen \fi
\vrule height \@halfwidth depth \@halfwidth width \@linelen
\ifnum \@xarg <0 \hskip -\@linelen \fi}

\def\@getlinechar(#1,#2){\@tempcnta#1\relax\multiply\@tempcnta 8
\advance\@tempcnta -9 \ifnum #2>0 \advance\@tempcnta #2\relax\else
\advance\@tempcnta -#2\relax\advance\@tempcnta 64 \fi
\char\@tempcnta}

\def\vector(#1,#2)#3{\@xarg #1\relax \@yarg #2\relax
\@linelen=#3\unitlength
\ifnum\@xarg =0 \@vvector
  \else \ifnum\@yarg =0 \@hvector \else \@svector\fi
\fi}

\def\@hvector{\@hline\hbox to 0pt{\@linefnt
\ifnum \@xarg <0 \@getlarrow(1,0)\hss\else
    \hss\@getrarrow(1,0)\fi}}

\def\@vvector{\ifnum \@yarg <0 \@downvector \else \@upvector \fi}

\def\@svector{\@sline
\@tempcnta\@yarg \ifnum\@tempcnta <0 \@tempcnta=-\@tempcnta\fi
\ifnum\@tempcnta <5
  \hskip -\wd\@linechar
  \@upordown\@clnht \hbox{\@linefnt  \if@negarg
  \@getlarrow(\@xarg,\@yyarg) \else \@getrarrow(\@xarg,\@yyarg) \fi}%
\else\@badlinearg\fi}

\def\@getlarrow(#1,#2){\ifnum #2 =\z@ \@tempcnta='33\else
\@tempcnta=#1\relax\multiply\@tempcnta \sixt@@n \advance\@tempcnta
-9 \@tempcntb=#2\relax\multiply\@tempcntb \tw@
\ifnum \@tempcntb >0 \advance\@tempcnta \@tempcntb\relax
\else\advance\@tempcnta -\@tempcntb\advance\@tempcnta 64
\fi\fi\char\@tempcnta}

\def\@getrarrow(#1,#2){\@tempcntb=#2\relax
\ifnum\@tempcntb < 0 \@tempcntb=-\@tempcntb\relax\fi
\ifcase \@tempcntb\relax \@tempcnta='55 \or
\ifnum #1<3 \@tempcnta=#1\relax\multiply\@tempcnta
24 \advance\@tempcnta -6 \else \ifnum #1=3 \@tempcnta=49
\else\@tempcnta=58 \fi\fi\or
\ifnum #1<3 \@tempcnta=#1\relax\multiply\@tempcnta
24 \advance\@tempcnta -3 \else \@tempcnta=51\fi\or
\@tempcnta=#1\relax\multiply\@tempcnta
\sixt@@n \advance\@tempcnta -\tw@ \else
\@tempcnta=#1\relax\multiply\@tempcnta
\sixt@@n \advance\@tempcnta 7 \fi\ifnum #2<0 \advance\@tempcnta 64 \fi
\char\@tempcnta}

\def\@vline{\ifnum \@yarg <0 \@downline \else \@upline\fi}

\def\@upline{\hbox to \z@{\hskip -\@halfwidth \vrule
  width \@wholewidth height \@linelen depth \z@\hss}}

\def\@downline{\hbox to \z@{\hskip -\@halfwidth \vrule
  width \@wholewidth height \z@ depth \@linelen \hss}}

\def\@upvector{\@upline\setbox\@tempboxa\hbox{\@linefnt\char'66}\raise
     \@linelen \hbox to\z@{\lower \ht\@tempboxa\box\@tempboxa\hss}}

\def\@downvector{\@downline\lower \@linelen
      \hbox to \z@{\@linefnt\char'77\hss}}

\thinlines

\newcount\@xarg
\newcount\@yarg
\newcount\@yyarg
\newcount\@multicnt
\newdimen\@xdim
\newdimen\@ydim
\newbox\@linechar
\newdimen\@linelen
\newdimen\@clnwd
\newdimen\@clnht
\newdimen\@dashdim
\newbox\@dashbox
\newcount\@dashcnt
 \catcode`@=12


\newbox\tbox
\newbox\tboxa

\def\leftzer#1{\setbox\tbox=\hbox to 0pt{#1\hss}%
     \ht\tbox=0pt \dp\tbox=0pt \box\tbox}

\def\rightzer#1{\setbox\tbox=\hbox to 0pt{\hss #1}%
     \ht\tbox=0pt \dp\tbox=0pt \box\tbox}

\def\centerzer#1{\setbox\tbox=\hbox to 0pt{\hss #1\hss}%
     \ht\tbox=0pt \dp\tbox=0pt \box\tbox}

%
\def\image(#1,#2)#3{\vbox to #1{\offinterlineskip
    \vss #3 \vskip #2}}


\def\leftput(#1,#2)#3{\setbox\tboxa=\hbox{%
    \kern #1\unitlength
    \raise #2\unitlength\hbox{\leftzer{#3}}}%
    \ht\tboxa=0pt \wd\tboxa=0pt \dp\tboxa=0pt\box\tboxa}

\def\rightput(#1,#2)#3{\setbox\tboxa=\hbox{%
    \kern #1\unitlength
    \raise #2\unitlength\hbox{\rightzer{#3}}}%
    \ht\tboxa=0pt \wd\tboxa=0pt \dp\tboxa=0pt\box\tboxa}

\def\centerput(#1,#2)#3{\setbox\tboxa=\hbox{%
    \kern #1\unitlength
    \raise #2\unitlength\hbox{\centerzer{#3}}}%
    \ht\tboxa=0pt \wd\tboxa=0pt \dp\tboxa=0pt\box\tboxa}

\unitlength=1mm
\input amssym.def
\input amssym.tex
\let\lgr=\longrightarrow
\def\lgrsim{\buildrel\sim\over\lgr}
\def\bpc#1{{\tenbf#1}\sevenbf}
\let\morinj=\hookrightarrow
\magnification=1200
\vsize=23true cm
\hsize=16true cm
\hoffset=0cm
\voffset=0cm

\centerline{{\bpc UN MOD\`ELE SEMI-STABLE DE LA VARI\'ET\'E DE} {\bpc SIEGEL DE GENRE }{\bf 3}}
\centerline{\bpc {}AVEC STRUCTURES DE NIVEAU DE TYPE ${\bf\Gamma_0}(p)$}
\vskip3mm
\centerline{par}
\vskip3mm
\centerline{Alain Genestier}
\vskip 2cm
{\bf 0\pointir  Introduction}
\vskip 5mm
 La vari\'et\'e 
de Siegel ${\cal S}(g,N,p)_{\bb C}={\cal H}_g/(\Gamma_0(p)\cap\Gamma (N))$ 
(o\`u ${\cal H}_g$ est le demi-espace de Siegel de genre $g$, $N$ est un entier
sup\'erieur ou \'egal \`a 3,
$p$ est un entier premier ne divisant pas $N$ et o\`u les sous groupes
$\Gamma_0(p)$ et
$\Gamma (N)$  du groupe  symplectique
${\rm Sp}(2g,{\bb Z})$  sont form\'es respectivement des matrices triangulaires 
inf\'erieures modulo $p$ et des matrices congrues \`a l'identit\'e modulo $N$) admet
un mod\`ele entier naturel ${\cal S}(g,N,p)$ sur ${\bb Z}[1/N]$ (cf. [dJ1] et aussi [CN] ou [R]
pour des cas analogues),  qui est lisse au-dessus de 
${\rm Spec}\,{\bb Z}[1/Np]$. 

Lorsque $g>1$, 
le mod\`ele ${\cal S}(g,N,p)$ n'est pas semi-stable en la place de mauvaise r\'eduction $p$.
Dans ce travail, nous allons cependant construire
pour de petites valeurs de $g$ ($g=2$ ou 3)
une "r\'esolution semi-stable" $\widetilde{\cal S}(g,N,p)$
de ${\cal S}(g,N,p)$, c'est \`a dire un morphisme propre 
$\widetilde{\cal S}(g,N,p)\lgr{\cal S}(g,N,p)$ 
qui sera un isomorphisme au-dessus de ${\rm Spec}\,{\bb Z}[1/Np]$
et dont la source sera semi-stable en $p$.

Lorsque $g$ vaut 2, la r\'esolution semi-stable $\widetilde{\cal S}(2,N,p)$ 
que nous obtiendrons a d\'ej\`a \'et\'e 
construite par de Jong  en utilisant une autre approche, et son r\'esultat
est  alors l\'eg\`erement plus pr\'ecis que le n\^otre (cf. [dJ2]
et la remarque 2 qui suit notre th\'eor\`eme 4.1). En revanche,
outre le r\'esultat d\'ej\`a obtenu pour $g\leq 3$,
la m\'ethode que nous utiliserons 
donne  l'espoir de construire  une r\'esolution semi-stable
de ${\cal S}(g,N,p)$ pour tout $g$ ; elle permet
aussi d'obtenir un mod\`ele semi-stable de certaines vari\'et\'es de Shimura
unitaires (cf. les remarques qui suivent les th\'eor\`emes 2.4.2 et 4.1). 

Nous allons maintenant donner une id\'ee de
la construction de notre r\'esolution semi-stable.

Le mod\`ele local ${\bb M}_g$ de Rapoport-Zink, qui est localement pour la topologie \'etale
isomorphe
\`a
${\cal S}(g,N,p)_{{\bb Z}_p}$, est d\'efini en termes d'alg\`ebre lin\'eaire (cf. le paragraphe
1) et est donc plus simple 
\`a utiliser que ${\cal S}(g,N,p)$. 

Nous construirons alors
(pour $g\leq 3$) une r\'esolution semi-stable $\widetilde{\cal L}_g\lgr {\bb M}_g$
du mod\`ele local ${\bb M}_g$. Plus pr\'ecis\'ement, le mod\`ele local s'envoie naturellement
vers la grassmannienne ${\cal L}$ des sous-modules lagrangiens du ${\bb Z}_p$-module 
symplectique "standard"  $V_0={\bb Z}_p$ (muni de la forme symplectique (1.2)) et
ce morphisme est birationnel (\`a un petit d\'etail pr\`es, cf. (2.3)).  
Nous construirons en fait (pour tout $g$) une suite d'\'eclatements 
$\widetilde{\cal L}_g\lgr {\cal L}_g$. Le th\'eor\`eme suivant (cf. th\'eor\`eme 2.4.2) 
fournira alors  la r\'esolution semi-stable d\'esir\'ee.
\th TH\'EOR\`EME
\enonce
Lorsque $g\leq 3$, le ${\bb Z}_p$-sch\'ema $\widetilde{\cal L}_g$ est semi-stable et 
l'application birationnelle $\widetilde{\cal L}_g\dashrightarrow{\bb M}_g$ est en fait un 
morphisme.
\endth
Le paragraphe 3 sera consacr\'e \`a la d\'emonstration de ce th\'eor\`eme.

Dans le dernier paragraphe, nous verrons comment la r\'esolution semi-stable 
$\widetilde{\cal L}_g\lgr{\bb M}_g$, qui est \'equivariante sous l'action du 
sch\'ema en groupes des sym\'etries du mod\`ele local, induit une r\'esolution semi-stable de 
${\cal S}(g,N,p)_{{\bb Z}_p}$.

Je remercie B. H. Gross qui m'a propos\'e de rechercher un mod\`ele semi-stable de 
${\cal S}(3,N,p)$ ; je le remercie aussi pour l'atmosph\`ere agr\'eable et stimulante 
que j'ai trouv\'ee \`a Harvard durant les trois mois qu'il m'a permis d'y passer. 
Mes remerciements vont aussi \`a A. J. de Jong, avec qui j'ai eu d'int\'eressantes 
conversations et qui m'a donn\'e acc\`es \`a certains de ses manuscrits non publi\'es, et \`a
M. Rapoport pour une enrichissante discussion que nous avons eue \`a Orsay, 
qui m'a confirm\'e l'int\'er\^et de la construction (2.4.1) (cf. 
les remarques sur les vari\'et\'es de Simura unitaires et leurs mod\`eles locaux 
qui suivent les deux th\'eor\`emes de cet article).
Je remercie finalement G. Laumon, qui m'a incit\'e \`a r\'ediger cet article.
 \vskip 5mm
{\bf 1\pointir  Rappels sur le mod\`ele local}
\vskip 2mm
1.1\pointir On va d'abord rappeler la d\'efinition du sch\'ema de modules des
vari\'et\'es ab\'eliennes principalement polaris\'ees de dimension $g$
avec structures de niveau de type $\Gamma(N)\cap\Gamma_0(p)$. 
On en donnera ensuite une formulation
\'equivalente plus commode.

Il r\'esulte de [GIT] qu'il existe un 
${\bb Z}[1/N]$-sch\'ema quasi-projectif ${\cal S}(g,N,p)$ dont l'ensemble des points \`a valeurs
dans tout  ${\bb Z}[1/N]$-sch\'ema noeth\'erien $S$ est l'ensemble des classes d'isomorphie de
quadruplets $(A,\lambda, \eta,H_\bullet)$, o\`u 

-- $A$ est un $S$-sch\'ema ab\'elien de 
dimension relative $g$ 

-- $\lambda$ est une polarisation principale de $A$

-- $\eta\ :\ ({\bb Z}/N{\bb Z})^g\oplus (\mu_N)^g\lgr A[N]$ est un isomorphisme
symplectique (on munit $({\bb Z}/N{\bb Z})^g\oplus (\mu_N)^g$ 
de la forme bilin\'eaire altern\'ee
\`a valeurs dans $\mu_N$
d\'efinie par la dualit\'e de Cartier entre les sch\'emas en groupes finis
${\bb Z}/N{\bb Z}$ et $\mu_N$, et $A[N]$ de l'accouplement de Weil d\'efini par $\lambda$).

-- $H_\bullet=(H_1\subset\cdots\subset H_g)$ est un drapeau de sous-sch\'emas en groupes 
finis localement libres de $A[p]$  tel que $H_i$ soit isotrope pour l'accouplement de 
Weil et d'ordre $p^i$.

Dans la suite de cet article, il  sera  plus commode d'utiliser 
la variante de pr\'esentation suivante
(cf. [dJ1], proposition 1.7) : l'ensemble 
${\cal S}(g,N,p)(S)$ est aussi celui des classes
d'isomorphie  de quadruplets
$(A_0\buildrel\alpha\over\rightarrow A_1\buildrel\alpha\over\rightarrow
\cdots\buildrel\alpha\over\rightarrow A_g,\lambda_0,\lambda_g,\eta)$, o\`u

-- $A_i$ est un $S$-sch\'ema ab\'elien de dimension relative $g$

-- $\alpha\ :\ A_i\lgr A_{i+1}$ est une isog\'enie de degr\'e $p$

-- $\lambda_0$ et $\lambda_g$ sont respectivement des polarisations principales de
$A_0$ et $A_g$ et rendent le diagramme
$$\diagram{A_0&\buildrel\alpha\over\rightarrow\cdots\buildrel\alpha\over\rightarrow& A_g\cr
\vfl{p\lambda_0}{}{4mm}&&\vfl{\lambda_g}{}{4mm}\cr
{}^t\!A_0&\buildrel{{}^t\!\alpha}\over\leftarrow
\cdots\buildrel{{}^t\!\alpha}\over\leftarrow&{}^t\!A_g
}$$
commutatif

-- $\eta\ :\ ({\bb Z}/N{\bb Z})^g\oplus (\mu_N)^g\lgr A_0[N]$ est un isomorphisme
symplectique

\noindent (on v\'erifie l'\'equivalence de ces deux d\'efinitions en prenant $A_0=A$ et $A_i=A/H_i$).
\vskip2mm 
1.2\pointir Le sous groupe d'Iwahori de ${\rm Sp}(2g,{\bb Q}_p)$ form\'e des matrices 
triangulaires inf\'erieures modulo $p$ est le groupe des points entiers d'un certain
${\bb Z}_p$-sch\'ema en groupes affine et lisse  ${\cal I}$, de fibre g\'en\'erique
${\rm Sp}(2g)_{{\bb Q}_p}$, dont on rappelera d'abord une construction.
On d\'efinira ensuite un $\cal I$-torseur 
${\cal T}(g,N,p)$ au-dessus de ${\cal S}(g,N,p)_{{\bb Z}_p}$.

Soient $V$ le ${\bb Q}_p$-espace vectoriel ${\bb Q}_p^{2g}$
et $(e_i)_{1\leq i\leq 2g}$ la base canonique de $V$. On munit 
$V$ de la forme 
altern\'ee $<\ ,>$ d\'efinie comme suit :
$$<x,y>={}^txJy$$
pour tout couple $x,y$ de vecteurs dans ${\bb Q}_p^{2g}$, o\`u 
$$J=\pmatrix{\hfill O&K\cr
-K&0}$$ 
et $K\in {\rm GL}(g,{\bb Z}_p)$ est la matrice   de la permutation $(g, g-1,\cdots,1)$
(les matrices $J$ et $K$ sont donc antidiagonales.

Soit d'autre part la suite de ${\bb Z}_p$-r\'eseaux $(V_i)_{1\leq i\leq g}$
d\'efinie par 
$$V_i=\Pi^{-i} {\bb Z}^{2g}_p$$
pour tout $i\in \bb Z$, o\`u $\Pi$ est la matrice
$$\pmatrix{
0&p\cr
{\rm Id}_{2g-1}&0
}$$
On note $\alpha$ : $V_i\morinj V_{i+1}$ l'inclusion \'evidente
et $\pi\ :\ V_i\lgrsim p^{-1}V_i=\kern-1mm =V_{i+2g}$ l'isomorphisme de p\'eriodicit\'e
(d\'efini par la multiplication par $p^{-1}$).
La forme altern\'ee $<\ ,>$ induit des formes altern\'ees non d\'eg\'en\'er\'ees
$$<\ ,>_{ig}\,=p^{i}<\ ,>\ :\ V_{ig}\times V_{ig}\lgr {\bb Z}_p\quad (i\in\{0,1\}).$$ 
Soit $\cal I$ le ${\bb Z}_p$-sch\'ema en groupes des automorphismes du syst\`eme 
$( V_0
\buildrel\alpha\over\rightarrow\cdots \buildrel\alpha\over\rightarrow V_i
\buildrel\alpha\over\rightarrow \cdots \buildrel\alpha\over\rightarrow V_g\,,
<\ ,\ >_0\,,<\  ,\ >_g)$. Il r\'esulte de ([dJ1], proposition 3.6) ou de ([RZ], 
theorem 3.16) que le sch\'ema en groupes $\cal I$ est lisse sur 
${\bb Z}_p$ ; ce sch\'ema en groupes est aussi celui d\'efini par la chambre
$(V_0\subset\cdots\subset V_i\subset\cdots\subset V_g)$ de l'immeuble de 
Bruhat-Tits de ${\rm Sp}(2g,{\bb Q}_p)$ (voir par exemple [T], 3.4 et 3.7).

Soient $(A_0\buildrel\alpha\over\rightarrow A_1\buildrel\alpha\over\rightarrow
\cdots\buildrel\alpha\over\rightarrow A_g,\lambda_0,\lambda_g,\eta)$ le quadruplet universel
au-dessus de ${\cal S}(g,N,p)$ et $\pi_i\ :\ A_i\lgr {\cal S}(g,N,p)$ le morphisme 
structural. Le premier faisceau de cohomologie de de Rham
$$M_i={\rm R}^1\pi_{(g-i)*}(\Omega^\cdot_{A_i/{\cal S}(g,N,p)})$$
est alors un ${\cal O}_{{\cal S}(g,N,p)}$-module localement 
libre de rang $2g$ ; les morphismes
$\alpha$ et les  polarisations $\lambda_\bullet$ induisent respectivement  des morphismes
$\alpha\ :\ M_i\lgr M_{i+1}$ et des formes symplectiques 
$<\ ,\ >_\bullet$ sur $M_0$ et $M_g$. Il r\'esulte de ([dJ1], proposition 3.6) 
(ou de [RZ], theorem 3.16) que le syst\`eme 
$((M_i)_{1\leq i\leq g},\alpha,<\ ,\ >_\bullet) $ est, localement
pour la topologie \'etale sur ${\cal S}(g,N,p)_{{\bb Z}_p}$, isomorphe  
\`a l'image r\'eciproque du syst\`eme 
$((V_i)_{1\leq i\leq g},\alpha,<\ ,\ >_\bullet) $. Le $\cal I$-torseur  ${\cal T}(g,N,p)$
au-dessus de ${\cal S}(g,N,p)_{{\bb Z}_p}$
est alors simplement celui des isomorphismes des deux syst\`emes :
$${\cal T}(g,N,p)={\rm Isom}\,(((M_i)_{1\leq i\leq g},\alpha, <\ ,\ >_\bullet),
((V_i)_{1\leq i\leq g},\alpha, <\ ,\ >_\bullet)_{{\cal S}(g,N,p)}).$$
\vskip2mm
1.3\pointir Le  ${\cal O}_{{\cal S}(g,N,p)}$-module localement 
libre de rang $g$ 
$$\omega_i=\pi_{(g-i)*}(\Omega^1_{A_i/{\cal S}(g,N,p)})$$
est localement un facteur direct de $M_i$. 
Les sous-modules $\omega_0$ et $\omega_g$ sont respectivement 
isotropes pour $<\ ,\ >_0$ et $<\ ,\ >_g$ ; on a de plus $\alpha(\omega_i)\subset
\omega_{i+1}=$.

Rapoport et Zink (cf. [R] ou [RZ]) introduisent alors le ${\bb Z}_p$-sch\'ema  
${\bb M}_g$ d\'efini par 
$$\displaylines{\hfill{\bb M}_g(R)=\{(\omega_i\subset 
V_i\otimes_{{\bb Z}_p}R)_{0\leq i\leq g}\mid
\hbox to 70mm{$\omega_i$ est un sous $R$-module localement\hfill}\cr
\hfill\hbox to 70mm{facteur direct de rang $g$ de  $V_i\otimes_{{\bb Z}_p}R$,\hfill}\cr
\hfill\hbox to 70mm{les sous-modules $\omega_0$ et $\omega_g$ sont isotropes\hfill}\cr
\hfill\hbox to 70mm{et on a $\alpha\otimes R\,(\omega_i) \subset\omega_{i+1}\}$ 
\hfill}}$$
pour toute ${\bb Z}_p$-alg\`ebre commutative $R$. La filtration de Hodge 
$(\omega_i\subset M_i)_{0\leq i\leq g}$ associ\'ee \`a l'objet universel
$$(A_0\buildrel\alpha\over\rightarrow A_1\buildrel\alpha\over\rightarrow
\cdots\buildrel\alpha\over\rightarrow A_g,\lambda_0,\lambda_g,\eta,\ 
\varphi_\bullet\ :\ (V_\bullet)_{{\cal T}(g,N,p)}\lgr M_\bullet)$$
sur ${\cal T}(g,N,p)$ induit un point $(\varphi_i^{-1}(\omega_i)\subset
(V_i)_{{\cal T}(g,N,p)})_{0\leq i\leq g}\in {\bb M}_g({\cal T}(g,N,p))$, 
et donc un morphisme $f\ :\ {\cal T}(g,N,p)\lgr {\bb M}_g$. 

Consid\'erons le diagramme 
$$\diagram{
&{\cal T}(g,N,p)
\cr\noalign{\vskip 5mm}
{\cal S}(g,N,p)_{{\bb Z}_p}\smash{\raise7mm\hbox{\kern-3mm $\scriptstyle  \rm pr$}}
&&{\bb M}_g
\smash{\raise7mm\hbox{\kern-8mm $\scriptstyle f$}} \cr
\fleche(20,11)\dir(-1,-1)\long{4}
\fleche(35,11)\dir(1,-1)\long{4}
}$$
Les ${\bb Z}_p$-sch\'emas ${\bb M}_g$ et ${\cal S}(g,N,p)_{{\bb Z}_p}$
sont tous deux de dimension relative $g(g+1)/2$.
Par ailleurs, il r\'esulte de la th\'eorie des d\'eformations des vari\'et\'es ab\'eliennes que 
le morphisme $f$ est lisse (cf. [Me], [RZ], 3. 29 ou [dJ1], proposition 4.5). A. J.  de Jong 
en tire alors le corollaire suivant (cf. [dJ1], corollary 4.7).
\th
PROPOSITION 1.3.1
\enonce
Pour tout point g\'eom\'etrique $s\ :\ {\rm Spec}\,k\lgr {\cal S}(g,N,p)_{{\bb Z}_p}$,
il existe un point g\'eom\'etrique $m\ :\ {\rm Spec}\,k\lgr {\bb M}_g$ et un isomorphisme
de ${\bb Z}_p$-sch\'emas
$V_s\lgrsim V_m$, o\`u $V_s$ (resp. $V_m$) est un voisinage \'etale de $s$ (resp. $m$)
dans ${\cal S}(g,N,p)$ (resp. ${\bb M}_g$). On peut de plus prendre 
$m=f(t)$, o\`u $t$ est n'importe quel point g\'eom\'etrique de ${\cal T}(g,N,p)$ au-dessus
de $s$. $\square$
\endth 
{\it Remarque} :  pour une autre d\'emonstration
de cet \'enonc\'e, voir aussi ([RZ], proposition 3. 33).

La proposition suivante implique alors que les sch\'emas ${\cal S}(g,N,p)$ et ${\bb M}_g$
sont localement ${\bb Z}_p$-isomorphes.
\th
PROPOSITION 1.3.2
\enonce
Le morphisme $f$ est surjectif.
\endth
{\it D\'emonstration} : Le morphisme $f$ est lisse et $\cal I$-\'equivariant, et son image est donc 
un ouvert $\cal I$-invariant de ${\bb M}_g$. Cet ouvert contient la fibre g\'en\'erique
de ${\bb M}_g$ (qui est form\'ee d'une seule orbite sous 
${\cal I}\otimes_{{\bb Z}_p}{\bb Q}_p$), et son compl\'ementaire est donc un ferm\'e de 
la fibre sp\'eciale de ${\bb M}_g$. Consid\'erons maintenant le point 
$m_0\in {\bb M}_g({\bb F}_p)$ d\'efini par 
$$\omega_i=[e'_{g+1},\cdots, e'_{2g}]\,,\ \forall 0\leq i\leq g$$
(o\`u l'on note $e'_j=\Pi^{-i}e_j\in V_j$ et $[e'_{g+1},\cdots, e'_{2g}]$ le sous ${\bb F}_p$
espace vectoriel de $V_i\otimes_{{\bb Z}_p}{\bb F}_p$ engendr\'e par $e'_{g+1},\cdots, e'_{2g}$).
La proposition (1.3.2) r\'esulte alors du lemme suivant.
\th 
LEMME 1.3.3
\enonce

{\rm 1}) Le point $m_0$ appartient \`a l'image du morphisme $f$.

{\rm 2}) Le point $m_0$ appartient \`a tout 
ferm\'e ${\cal I}\otimes_{{\bb Z}_p}{\bb F}_p$-invariant de la fibre sp\'eciale de ${\bb M}_g$.
\endth
{\it D\'emonstration du lemme} : 1) (cf. [dJ1,\S 5]). Soit $E$ une courbe elliptique supersinguli\`ere
sur $\bar{\bb F}_p$, munie de $\alpha_p\morinj E$ et prenons $A_i=(E/\alpha_p)^i\times
E^{g-i}$. On d\'efinit les morphismes $\alpha\ :\ A_i\lgr A_{i+1}$ de mani\`ere \'evidente
et on munit $A_0$ et $A_g$ de la polarisation produit. Il est facile de v\'erifier
que ce\c ci d\'efinit un point $s_0$ de ${\cal S}(g,N,p)(\bar{\bb F}_p)$ et que pour tout
$t_0\in{\cal S}(g,N,p)(\bar{\bb F}_p)$ au-dessus de $s_0$, on a $f(t_0)=m_0$ 
(on consid\`ere i\c ci $m_0$ comme un point de ${\bb M}_g(\bar{\bb F}_p)$ ; le point
$m_0$ est invariant sous l'action de ${\cal I}\otimes_{{\bb Z}_p}\bar{\bb F}_p$ 
de sorte qu'il suffit de choisir un point $t_0$ au-dessus de $s_0$).

2) Soit ${\bb V}_i$ le ${\bb F}_p[t]$-module $T^{-i}\,{\bb F}_p[t]^{2g}$, o\`u
l'on note
$$T=\pmatrix{0&t\cr
{\rm Id}_{2g-1}&0}.$$
On d\'efinit des morphismes $\alpha\ :\ {\bb V}_i\morinj{\bb V}_{i+1}$,
$\pi\ :\ {\bb V}_i\lgrsim t^{-1}{\bb V}_{i+2g}$ et des formes altern\'ees
$<\ ,\ >_{ig}\ ;\ {\bb V}_{ig}\times{\bb V}_{ig}\lgr {\bb F}_p[t]$ de mani\`ere analogue \`a  
ceux de (1.2). On a une identification naturelle
${\bb V}_i/t{\bb V}_i\lgrsim V_i\otimes_{{\bb Z}_p}{\bb F}_p$ compatible aux morphismes 
$\alpha$, $\pi$, et aux formes altern\'ees $<\ ,\ >_{ig}$.
Consid\'erons alors la suite de ${\cal O}_{{\bb M}_g}[t]$-r\'eseaux 
$(\Omega_i\subset{\cal O}_{{\bb M}_g}[t,t^{-1}]^{2g})_{0\leq i\leq g}$ d\'efinie par 
$\Omega_i=T^{-g}({\bb V}_i\otimes_{{\bb F}_p}{\cal O}_{{\bb M}_g}
\lgr  V_i\otimes_{{\bb Z}_p}{\cal O}_{{\bb M}_g})^{-1}(\omega_i)$. Ceci d\'efinit
un plongement de ${\bb M}_g\otimes_{{\bb Z}_p}{\bb F}_p$ 
dans la vari\'et\'e de drapeaux affine $\widetilde{\cal B}$ (cf. [KL,\S 5]) de
${\rm Sp}(2g)_{{\bb F}_p}$ (la fibre sp\'eciale de ${\bb M}_g$ 
est alors une g\'en\'eralisation naturelle
des sous-sch\'emas $X_i$ de [KL, 5.2] : en remarquant que $T^g$ est une racine carr\'ee de 
$t.{\rm Id}_{2g}$, la fibre sp\'eciale de ${\bb M}_g$ est alors le sch\'ema $X_{1/2}$).
Le groupe ${\cal I}\otimes_{{\bb Z}_p}{\bb F}_p$ est naturellement un quotient
de l'ind-groupe alg\'ebrique not\'e $\widetilde B$ dans [KL] et le plongement
${\bb M}_g\morinj \widetilde {\cal B}$ que nous venons de d\'ecrire est alors 
$\widetilde{ B}$-\'equivariant. Son image est donc une r\'eunion finie de 
sous-vari\'et\'es de Schubert g\'en\'eralis\'ees, et la deuxi\`eme partie du lemme r\'esulte
alors  de [KL, 5.2] (la $\widetilde B$-orbite $\{m_0\}\subset\widetilde{\cal B}$, qui 
correspond \`a l'unit\'e du groupe de Weyl affine de ${\rm Sp}(2g)$, est incluse 
dans toutes les sous-vari\'et\'es de Schubert). Ceci termine la d\'emonstration du lemme 1.3.3,
et donc aussi celle de la proposition 1.3.2. $\square$

De m\^eme que Rapoport et Zink (cf. [R] ou [RZ]), nous appelerons mod\`ele local le sch\'ema 
${\bb M}_g$.
\vfill\eject
{\bf 2\pointir  Une r\'esolution semi-stable du mod\`ele local (pour $g$= 2 ou 3)}
\vskip2mm
2.1\pointir On note $G$ le ${\bb Z}_p$-sch\'ema en groupes 
${\rm Sp}(V_0,<\ ,\ >_0)$ ; on note respectivement $P$ et $B$ les paraboliques qui 
stabilisent le sous-module isotrope $[e_{g+1},\cdots,e_{2g}]$ et le drapeau isotrope
$[e_{2g}]\subset[e_{2g-1},e_{2g}]\subset \cdots\subset [e_{g+1},\cdots , e_{2g}]$ 
(lorsque $S$ est
une partie de $V={\bb Q}_p^{2g}$, $[S]$ d\'esigne le ${\bb Z}_p$-module engendr\'e par $S$).
Le sch\'ema en groupes $B$ n'est autre que celui des matrices triangulaires inf\'erieures 
de ${\rm Sp}(2g)_{{\bb Z}_p}$ et $P$ est le parabolique  de Siegel \'evident
contenant ce sous-groupe de Borel.

On note ${\cal L}_g$ la grassmannienne des sous-modules lagrangiens de $V_0$, d\'efinie par
$$\displaylines{\qquad{\cal L}_g(R)=\{\omega_0\subset V_0\otimes_{{\bb Z}_p}R\mid
\hbox{$\omega_0$ est un sous $R$-module isotrope}\hfill\cr
\hfill\hbox{et localement facteur direct de rang 
$g$ de $V_0\otimes_{{\bb Z}_p}R$}\}\!\!\!\qquad}$$
pour toute ${\bb Z}_p$-alg\`ebre commutative $R$. 
C'est un sch\'ema lisse de dimension relative
$g(g+1)/2$ sur ${\bb Z}_p$. La grassmannienne ${\cal L}_g$ est munie d'une action 
localement transitive du sch\'ema en groupes $G$ ; le stabilisateur du point
$\omega_0=[e_{g+1}\cdots e_{2g}]$ n'est autre que $P$, 
de sorte que ${\cal L}_g$ peut aussi \^etre vue
comme le quotient $G/P$. 
\vskip2mm
2.2\pointir Nous allons maintenant faire quelques rappels (inspir\'es de [LS] et de [MS])
sur les cellules de Schubert
et les sous-vari\'et\'es de Schubert de la grassmannienne ${\cal L}_g$.
 
Soit $D_\bullet =(0)=D_{2g}\subset D_{2g-1}\subset\cdots\subset D_0=V_0$ le drapeau
d\'efini par $D_i=[e_{i+1},e_{i+2},\cdots,e_{2g}]$ 
(de sorte que $D_\bullet$ est fixe par $B$ et que
${\rm rg}_{{\bb Z}_p}(D_i)=2g-i$). 
Lorsque $S=\{\lambda_1,\cdots,\lambda_g\}$ (o\`u l'on suppose que
$\lambda_1<\lambda_2<\cdots<\lambda_g$)
est une partie \`a $g$ \'el\'ements de $\{1,\cdots,2g\}$, soit ${\cal C}_S$ 
le sous-sch\'ema localement ferm\'e de ${\cal L}_g$ d\'efini par
$$\displaylines{\qquad{\cal C}_S(R)=\{(\omega_0\subset V_0)\in {\cal L}_g(R)\mid
\hbox{localement sur ${\rm Spec}\, R$ pour la topologie}\hfill\cr
\hphantom{\qquad{\cal C}_S(R)=\{(\omega_0\subset V_0)\in {\cal L}_g(R)\mid}
\ \hbox{de Zariski, il existe une base $(s_i)_{1\leq i\leq g}$ }\hfill\cr
\hphantom{\qquad{\cal C}_S(R)=\{(\omega_0\subset V_0)\in {\cal L}_g(R)\mid}
\ \hbox{de $\omega_0$ telle que $s_i$ soit une section (locale) }\hfill\cr
\hphantom{\qquad{\cal C}_S(R)=\{(\omega_0\subset V_0)\in {\cal L}_g(R)\mid}
\ \hbox{de $D_{i-1}$ engendrant le module inversible $D_{i-1}/D_i$}\}\hfill}$$
pour toute ${\bb Z}_p$-alg\`ebre commutative $R$ ; on dit que ${\cal C}_S$
est la cellule de Schubert associ\'ee \`a $S$ (cf. par exemple [LS]). 
Le sch\'ema ${\cal C}_S$ est non vide
si et seulement si $S$ est une partie  {\it  totalement isotrope}
(c'est \`a dire telle que $[e_i\,,i\in S]$ soit totalement isotrope)
de $\{1,\cdots,2g\}$ ; le sch\'ema ${\cal C}_S$ est alors la $B$-orbite du point
de ${\cal L}_g$ d\'efini par le sous-module $[e_i\,,i\in S]$ de $V_0$. 
On d\'efinit la vari\'et\'e de Schubert ${\cal L}_S$ associ\'ee \`a $S$  comme 
l'adh\'erence sch\'ematique de 
la cellule de Schubert ${\cal C}_S$ (cf. par exemple [LS]).

La proposition suivante r\'esume quelques propri\'et\'es des cellules et des vari\'et\'es de
Schubert (pour une d\'emonstration, cf.par exemple [LS]).
\th
PROPOSITION 2.2.1
\enonce

{\rm 1}) Toutes les orbites de ${\cal L}_g$ sous l'action de $B$ 
sont en fait de la forme ${\cal C}_S$,
et on a donc une stratification localement ferm\'ee $B$-\'equivariante
$${\cal L}_g=\coprod_S{\cal C}_S$$
o\`u $S$ parcourt l'ensemble $\goth L$ des parties totalement isotropes de 
$\{1,\cdots,2g\}$. 

{\rm 2}) La vari\'et\'e de Schubert ${\cal L}_S$ a alors
une stratification localement ferm\'ee induite de la forme
${\cal L}_S=\coprod_{S'}{\cal C}_{S'}$, o\`u $S'$ d\'ecrit une partie ${\goth L}_S$ de $\goth L$.
En particulier, on a ${\cal L}_{S'}\subset {\cal L}_{S}$ si et seulement si
$S'\in {\goth L}_S$.

{\rm 3}) La partie ${\goth L}_S$ est en fait constitu\'ee des parties totalement isotropes 
$S'=\{\lambda'_1,\cdots,\lambda'_g\}\in {\goth L}$ 
(dont on a rang\'e les \'el\'ements $\lambda'_i$ dans l'ordre croissant)
telles que 
$\lambda_i\leq\lambda'_i\,,\ \forall i$. 

{\rm 4}) Le ${\bb Z}_p$-sch\'ema ${\cal C}_S$ est isomorphe \`a un espace
affine ${\bb A}_{{\bb Z}_p}^{\ell(S)}$. Cet espace affine est de dimension
$$\ell(S)=r(g+1)-\sum_{i=1}^r\lambda_i\,,$$
o\`u l'on note $r$ le nombre d'\'el\'ements de $S\cap \{1,\cdots,g\}$.
$\square$\endth

{\it Exemples} : on a  pour $g=2$ et $g=3$ les deux diagrammes suivants
$$\diagram{{\cal L}_{\{3,4\}}\hfl{}{}{6mm}{\cal L}_{\{2,4\}}
\hfl{}{}{6mm}{\cal L}_{\{1,3\}}\hfl{}{}{6mm}{\cal L}_{\{1,2\}}
={\cal L}_2
}$$
$$\diagram{&\kern-6mm{\cal L}_{\{1,4,5\}}&\cr\noalign{\vskip5mm}
{\cal L}_{\{4,5,6\}}\hfl{}{}{6mm}{\cal L}_{\{3,5,6\}}
\hfl{}{}{6mm}{\cal L}_{\{2,4,6\}}
&&\kern-2mm
{\cal L}_{\{1,3,5\}}\hfl{}{}{6mm}{\cal L}_{\{1,2,4\}}\hfl{}{}{6mm}
{\cal L}_{\{1,2,3\}}={\cal L}_3\cr\noalign{\vskip5mm}
&\kern-6mm{\cal L}_{\{2,3,6\}}&\cr
\fleche(39,18)\dir(1,1)\long{4}
\fleche(39,11)\dir(1,-1)\long{4}
\fleche(51,7)\dir(1,1)\long{4}
\fleche(51,22)\dir(1,-1)\long{4}}$$ 
o\`u les fl\`eches figurent des immersions ferm\'ees de codimension $1$
(de sorte que la dimension relative de ${\cal L}_S$ sur ${\bb Z}_p$
n'est autre que la distance de ${\cal L}_S$ au sommet le plus \`a gauche
du graphe). En fait, le premier de ces deux diagrammes n'est autre que le sous-diagramme
$$\diagram{{\cal L}_{\{4,5,6\}}\hfl{}{}{6mm}{\cal L}_{\{3,5,6\}}
\hfl{}{}{6mm}{\cal L}_{\{2,4,6\}}\hfl{}{}{6mm}{\cal L}_{\{2,3,6\}}
}$$
du second.
\vskip2mm
2.3\pointir On note ${\rm Gr}_i$ la grassmannienne des sous-modules 
(localement facteurs directs) de rang $g$ de $V_i$. On consid\`ere le morphisme
$$\delta\ :\ {\cal L}_g\otimes_{{\bb Z}_p}{\bb Q}_p\lgr
\prod_{0\leq i \leq g}{\rm Gr}_i\otimes_{{\bb Z}_p}{\bb Q}_p$$
d\'efini par 
$$\delta (\omega_0)=(\omega_0\subset
V\otimes_{{\bb Q}_p}R=\!=V_i\otimes_{{\bb Z}_p}R)_{0\leq i\leq g}$$
pour toute ${\bb Q}_p$-alg\`ebre commutative  R. L'image de $\delta$ est un sous-sch\'ema
localement ferm\'e de ${\bb M}_g\subset\prod_{0\leq i \leq g}{\rm Gr}_i$ 
et on notera ${\cal M}_g$
son adh\'erence de Zariski. L'action naturelle du sch\'ema en groupes $\cal I$ sur
${\bb M}_g$ stabilise ${\cal M}_g$ et d\'efinit donc une action de
$\cal I$ sur ${\cal M}_g$. La projection sur le facteur ${\rm Gr}_0$ d\'efinit
un morphisme $\cal I$-\'equivariant ${\rm pr}_0\ :\ {\cal M}_g \lgr {\cal L}_g$
(on fait agir $\cal I$ sur ${\cal L}_g$ via le morphisme ${\cal I}\lgr G$ \'evident.
\vskip2mm
{\it Remarques} : 1) La fibre g\'en\'erique de l'immersion ferm\'ee tautologique
${\cal M}_g\morinj {\bb M}_g$ est un isomorphisme.
Le sch\'ema ${\cal M}_g$ est donc la "platification" du sch\'ema ${\bb M}_g$ (i.e. le sch\'ema
obtenu \`a partir de ${\bb M}_g$ en supprimant la $p$-torsion ; cette op\'eration
est un cas (tr\`es) particulier de [GR, 5.2]).
En fait, Rapoport et Zink  conjecturent que le sch\'ema ${\bb M}_g$ est d\'ej\`a plat
(cf. les trois lignes qui pr\'ec\`edent le paragraphe 3.36 de [RZ]).

2) L'op\'eration de platification de 
la remarque pr\'ec\'edente est locale pour la topologie \'etale, et ${\cal M}_g$ est donc
aussi un mod\`ele local pour le platifi\'e de ${\cal S}(g,N,p)$.
\vskip2mm
2.4\pointir  Lorsque $g\geq 2$, le sch\'ema ${\cal M}_g$ n'est pas semi-stable. 
En fait, la fibre sp\'eciale de 
${\cal M}_g$ n'est pas un diviseur \`a croisements normaux et lorsque $g\geq 3$,
${\cal M}_g$ n'est m\^eme pas localement d'intersection compl\`ete.

On aimerait disposer d'une
"r\'esolution semi-stable $\cal I$-\'equivariante
$\widetilde{\cal M}_g\lgr {\cal M}_g$ de ${\cal M}_g$", c'est \`a dire d'un morphisme de 
${\bb Z}_p$-sch\'emas $\widetilde{\cal M}_g\lgr {\cal M}_g\,$, 
$\cal I$-\'equivariant et propre,
dont la fibre g\'en\'erique est un isomorphisme et dont la source $\widetilde{\cal M}_g$
est semi-stable
(par quoi l'on entend que, localement pour la topologie \'etale, celle \c ci est 
${\bb Z}_p$-isomorphe \`a ${\rm Spec\,}{\bb Z}_p[t_1,\cdots,t_n]/(t_1\cdots t_r-p)$).
Pour $g=$ 2 ou 3, nous allons effectivement en  construire une.

Remarquons tout d'abord que par composition avec le morphisme $\cal I$-\'equivariant
${\rm pr}_0\ :\ {\cal M}_g\lgr{\cal L}_g$ (2.3), la r\'esolution semi-stable
$\cal I$-\'equivariante
$\widetilde{\cal M}_g\lgr {\cal M}_g$ d\'esir\'ee va induire
un morphisme propre $\cal I$-\'equivariant $\widetilde{\cal M}_g\lgr {\cal L}_g$,
dont la fibre g\'en\'erique sera un isomorphisme. Ceci sugg\`ere de construire 
$\widetilde{\cal M}_g$ \`a partir de ${\cal L}_g$ \`a l'aide d'une suite
d'\'eclatements $\cal I$ \'equivariants. Nous allons donc construire, pour tout $g$,
une suite d'\'eclatements $\cal I$-\'equivariants $\widetilde{\cal L}_g\lgr {\cal L}_g$.
Nous verrons ensuite que pour $g=$ 2 ou 3, le candidat  $\widetilde{\cal L}_g$
fournit une r\'esolution semi-stable $\cal I$-\'equivariante 
$\widetilde{\cal L}_g\lgr {\cal M}_g$ et m\'erite donc le nom de $\widetilde{\cal M}_g$.
\vskip2mm
Construction 2.4.1\pointir On note  
$\overline{\cal L}_S={\cal L}_S\otimes_{{\bb Z}_p}{\bb F}_p$ la fibre 
sp\'eciale de ${\cal L}_S$. Lorsque $0\leq i \leq g(g+1)/2$, on note 
${\goth L}_i=\{S\in {\goth L}\mid l(S)=i\}$ l'ensemble des $S$ tels que 
la sous-vari\'et\'e de Schubert $\overline{\cal L}_S\subset\overline{\cal L}_g$ 
soit de dimension relative $i$ sur ${\bb F}_p$. La construction de 
$\widetilde{\cal M}_g$ va alors se faire en plusieurs \'etapes.

$\underline{\hbox{Etape 1}}$ : on \'eclate ${\cal L}_g$ le long de l'unique sous-vari\'et\'e 
de Schubert de dimension nulle, $\overline{\cal L}_{\{g+1,\cdots,2g\}}$,  
de sa fibre sp\'eciale. On obtient ainsi un ${\bb Z}_p$-sch\'ema 
$\widetilde{\cal L}_{\leq 0}$

$\underline{\hbox{Etape 2}}$ : on \'eclate $\widetilde{\cal L}_{\leq 0}$ le long
du transform\'e strict de l'unique sous-vari\'et\'e 
de Schubert de dimension 1, $\overline{\cal L}_{\{g,g+2,\cdots,2g\}}$, 
de la fibre sp\'eciale de ${\cal L}_g$. On obtient ainsi un ${\bb Z}_p$-sch\'ema 
$\widetilde{\cal L}_{\leq 1}$

...

$\underline{\hbox{Etape $i$}}$ : lorsque $S\in {\goth L}_{i-1}$, notons 
$(\overline{\cal L}_S)\,\widetilde{}\kern-3pt{}_{\leq i-2}$ le transform\'e strict de 
$\overline{\cal L}_S$ par la cha\^ine d'\'eclatements 
$\widetilde{\cal L}_{i-2}\lgr\cdots\lgr\widetilde{\cal L}_{\leq 0}\lgr{\cal L}_{g}$.
On \'eclate  $\widetilde{\cal L}_{i-2}$ le long de 
$\bigcup_{S\in {\goth L}_{i-1}}\,(\overline{\cal L}_S)\,\widetilde{}\kern-3pt{}_{\leq i-2}$.
On obtient ainsi 
un ${\bb Z}_p$-sch\'ema $\widetilde{\cal L}_{\leq i-1}$

...

$\underline{\hbox{Etape $g(g+1)/2$}}$ : on \'eclate 
$\widetilde{\cal L}_{\leq g(g+1)/2-2}$ le long
du transform\'e strict (par les \'eclatements pr\'ec\'edents)
de l'unique sous-vari\'et\'e 
de Schubert de dimension $g(g+1)/2-1$,
$\overline{\cal L}_{\{1,\cdots,g-1,g+1\}}$,
de la fibre sp\'eciale de ${\cal L}_g$. On obtient ainsi le ${\bb Z}_p$-sch\'ema 
$\widetilde{\cal L}_g=\widetilde{\cal L}_{\leq g(g+1)/2-1}$.  

Remarquons que les sous-sch\'emas $\overline{\cal L}_S$ de ${\cal L}_g$ sont 
invariants sous l'action de $\cal I$ (ils sont invariants sous l'action de
$B_{{\bb F}_p}\subset G_{{\bb F}_p}$ sur $\overline{\cal L}_g$).
Le sch\'ema $\widetilde{\cal L}_g$ qu'on vient de construire est donc muni
d'une action naturelle de $\cal I$, pour laquelle le morphisme \'evident
${\rm pr}\ :\ \widetilde{\cal L}_g\lgr {\cal L}_g$ est $\cal I$-\'equivariant.
Ceci termine la construction.
\vskip2mm
{\it Remarque} : On d\'emontre que lorsque $g\leq 3$, la r\'eunion 
$\bigcup_{S\in {\goth L}_{i-1}}\,(\overline{\cal L}_S)\,\widetilde{}\kern-3pt{}_{\leq i-2}$
est en fait une union disjointe de vari\'et\'es lisses sur ${\bb F}_p$. Les \'eclatements qui
interviennent dans la construction ont alors tous des centres lisses.
\vskip2mm
Le th\'eor\`eme suivant a d\'ej\`a \'et\'e annonc\'e dans l'introduction.
\th TH\'EOR\`EME 2.4.2
\enonce Supposons que l'entier $g$ soit inf\'erieur ou \'egal \`a $3$.
Les propri\'et\'es suivantes sont alors v\'erifi\'ees.

1) Le ${\bb Z}_p$-sch\'ema $\widetilde{\cal L}_g$ est semi-stable.

2) Il existe un (et un seul) morphisme 
$R\ :\ \widetilde{\cal L}_g\lgr{\cal M}_g$ rendant le diagramme
$$\kern25mm\diagram{
\widetilde{\cal L}_g & {\hfl{R}{}{14mm}} &{\cal M}_g
\cr\noalign{\vskip2mm}
\vfl{\rm pr}{}{10mm}\cr\noalign{\vskip2mm}
{\cal L}_g&&&\smash{\raise14mm\hbox to 0mm{\kern-15mm $\scriptstyle {\rm pr}_0 
\ (cf.\ 2.3)$\hss}}
\cr
\fleche(25,25)\dir(-1,-1)\long{18}
}$$
commutatif. Ce morphisme $R$ est de plus propre et $\cal I$-\'equivariant. $\square$
\endth
Lorsque $g\leq 3$, nous noterons donc encore $\widetilde{\cal M}_g$ le
sch\'ema $\widetilde{\cal L}_g$.
\vskip2mm
{\it Remarques} : 0) Le cas o\`u $g=1$ 
que nous n'avons pas mentionn\'e jusqu'i\c ci
est trivial et la r\'esolution semi-stable obtenue est un isomorphisme.

1) Il semble naturel de se demander si l'\'enonc\'e 2.4.2 reste
vrai lorsqu'on supprime la restriction $g\leq 3$.

 2) En utilisant le fait que ${\cal L}_2$ 
n'est autre que la sous-vari\'et\'e de Schubert ${\cal L}_{\{2,3,6\}}$ de 
${\cal L}_3$, on voit qu'une preuve du th\'eor\`eme (2.4.2) dans le cas o\`u $g=2$ est
en fait "contenue" dans la preuve que nous donnerons pour $g=3$.

La r\'esolution semi-stable $\widetilde{\cal M}_2\lgr{\cal M}_2$ a d'autre part d\'ej\`a
\'et\'e obtenue par de Jong (cf. [dJ2])
par une autre m\'ethode. Plus pr\'ecis\'ement, de Jong \'eclate ${\cal M}_2$ le long
d'une composante irr\'eductible de sa fibre sp\'eciale (celle not\'ee $Z_{00}$
dans [dJ1, \S 5] ; on peut aussi choisir $Z_{11}$). 
On peut v\'erifier que la r\'esolution semi-stable ainsi obtenue
co\"incide avec la n\^otre.

3) Soit $W_0\buildrel\alpha\over\morinj W_1\buildrel\alpha\over\morinj\cdots
\buildrel\alpha\over\morinj W_{n-1}\buildrel\alpha\over\morinj W_n=p^{-1}W_0$
une cha\^ine de ${\bb Z}_p$-r\'eseaux dans $W={\bb Q}^n_p$. Pour $1\leq r\leq n-1$, 
consid\'erons le sch\'ema
$${\bb M}_{U(r,n-r)}=\{(\omega_i\subset W_i)_i\in \prod_{0\leq i\leq n-1}
{\rm Gr}(r,W_i)\mid \alpha(\omega_i)\subset \omega_{i+1},\forall i\}$$
(o\`u on note ${\rm Gr}(r,W_i)$ la grassmannienne des sous-modules localement 
facteurs directs de rang $r$ de $W_i$), qui est lui aussi un mod\`ele local
d'une vari\'et\'e de Shimura (pour un groupe de similitudes unitaires,  de type $GU(r,n-r)$ en la  
place archim\'edienne ,  d\'efini par une 
extension quadratique d\'eploy\'ee
en $p$, et pour des structures de niveau  
convenables ; pour plus de pr\'ecisions cf. [R, \S 3]).

La construction (2.4.1) admet dans ce cadre un analogue \'evident. 
L'analogue du th\'eor\`eme 
(2.4.2) est d\'ej\`a connu lorsque $r=1$ (ou $r=n-1$)(le morphisme $R$ est alors un 
isomorphisme : cf. dans [BC] ou [R] 
l'expos\'e de la construction de 
Deligne pour le sch\'ema formel $\widehat\Omega^n$).
On peut \'egalement v\'erifier cet analogue
dans le cas particulier o\`u $r=2$ et $n=4$ (la d\'emonstration est 
semblable \`a celle qui figure dans le paragraphe suivant). 

\vskip 5mm 
{\bf 3\pointir D\'emonstration du th\'eor\`eme 2.4.2}
\vskip 2mm
Nous nous contenterons de d\'emontrer ce th\'eor\`eme dans
le cas particulier o\`u $g$ vaut 3, qui est le plus difficile. Nous omettrons alors l'indice 3 de 
${\cal L}_3,\ {\cal M}_3,\ \cdots$
\vskip2mm
3.1\pointir Soit $B_-$ le sous-sch\'ema en groupes de G form\'e des matrices triangulaires 
sup\'erieures
($B_-$ est donc le le sch\'ema en groupes de Borel oppos\'e \`a $B$).
Soit $x_0=([e_4,e_5,e_6]\subset V_0)\in {\cal L}({\bb Z}_p)$
l'unique point de ${\cal L}_{\{4,5,6\}}$. La grosse cellule 
de Schubert oppos\'ee ${\cal L}^0$ (cf. [MS]) est la $B_-$-orbite $B_-x_0\subset \cal L$.
C'est un sous-sch\'ema ouvert de $\cal L$, qui s'identifie au translat\'e de
la cellule de Schubert ouverte ${\cal C}_{\{1,2,3\}}$ par l'\'el\'ement le plus
long $J$ (cf. (1.2)) de $W$ et est 
en particulier isomorphe \`a l'espace affine
${\bb A}^6_{{\bb Z}_p}$. Un isomorphisme "naturel"
${\bb A}^6_{{\bb Z}_p}\lgrsim{\cal L}^0$ s'obtient de la mani\`ere suivante.
Consid\'erons les points de ${\bb A}^6_{{\bb Z}_p}$ comme des 
matrices sym\'etriques $A=(a^j_i)_{1\leq i,j\leq 3}\,$. La matrice
$\pmatrix{A\cr
K\cr}$ (cf. (1.2))  est alors de rang 3, et on v\'erifie ais\'ement que
son image d\'efinit un point de ${\cal L}^0$. 

L'ouvert ${\cal L}^0$ contient le point $x_0$, qui est adh\'erent 
\`a toutes les cellules de Schubert ${\cal C}_S$. L' intersection
${\cal C}_S\cap {\cal L}^0$ est donc un ouvert non vide de ${\cal C}_S$ ;
en particulier ${\cal L}^0$ rencontre toutes les orbites de $B$. Le morphisme
$B\times_{{\bb Z}_p}{\cal L}^0\lgr\cal L$ d\'efini par l'action (2.2) de 
$B\subset G$ sur $\cal L$ est donc surjectif. A fortiori, le morphisme
${\cal I}\subset G$ sur $\cal L$ est alors surjectif ; nous dirons que l'ouvert 
${\cal L}^0$ est $\cal I$-{\it saturant}. 

Pour v\'erifier la premi\`ere partie du th\'eor\`eme (2.4.2) 
(la semi-stabilit\'e de $\widetilde{\cal M}$), il suffit manifestement de le faire 
au-dessus d'un ouvert $\cal I$-saturant de $\cal L$ (ou de $\widetilde{\cal L}_{\leq i}$).
Le lemme suivant nous permettra de r\'eduire de m\^eme  la d\'emonstration de la deuxi\`eme
partie de ce th\'eor\`eme \`a celle d'un \'enonc\'e concernant la restriction de 
$\widetilde{\cal M}$ au-dessus d'un ouvert saturant de $\widetilde{\cal L}_{\leq i}$.

\th LEMME 3.1.1
\enonce Soient $\Gamma$ un ${\bb Z}_p$-sch\'ema en groupes lisse et $X,Y,Z$ trois 
${\bb Z}_p$-sch\'emas munis d'une action de $\Gamma$. Soient
$f\, :\, X\lgr Z$ et $g\, :\,Y\lgr Z$ deux morphismes $\Gamma$-\'equivariants.
 Supposons que le sch\'ema $X$ est r\'eduit,  
que le morphisme $g$ est s\'epar\'e et qu'il existe un ouvert dense  de $Z$,
d'image r\'eciproque dense dans $X$,
au-dessus duquel $g$ est une immersion. 

Soient $U$ un ouvert $\Gamma$-saturant  de $X$ et supposons qu'il existe un morphisme
 $h_U\,:\,U\lgr Y$  rendant le diagramme
$$\kern25mm\diagram{
{U} & {\hfl{h_U}{}{16mm}} &{Y}
\cr\noalign{\vskip2mm}
\vfl{f}{}{10mm}\cr\noalign{\vskip2mm}
{Z}&&&\smash{\raise14mm\hbox to 0mm{\kern-15mm $\scriptstyle g 
$\hss}}
\cr
\fleche(25,25)\dir(-1,-1)\long{18}
}$$
commutatif (sous les hypoth\`eses pr\'ec\'edentes, il existe au plus un tel morphisme).

Il existe alors un (et un seul) morphisme 
$h\,:\,X\lgr Y$ rendant le diagramme 
$$\kern25mm\diagram{
{X} & {\hfl{h}{}{16mm}} &{Y}
\cr\noalign{\vskip2mm}
\vfl{f}{}{10mm}\cr\noalign{\vskip2mm}
{Z}&&&\smash{\raise14mm\hbox to 0mm{\kern-15mm $\scriptstyle g 
$\hss}}
\cr
\fleche(25,25)\dir(-1,-1)\long{18}
}$$
commutatif. Le morphisme $h$ est $\Gamma$-\'equivariant, et il est propre lorsque $f$ 
est propre.

\endth
{\it D\'emonstration} : Notons $v\ :\ V=\Gamma\times_{{\bb Z}_p}U
\lgr X$ le morphisme d\'efini par l'action de $\Gamma$ sur $X$ 
et consid\'erons le morphisme 
$\ h_V\,:\,V\lgr Y$ d\'efini par 
$(\gamma,u)\mapsto \gamma.  h_U(u)$.
Le morphisme $\,v$ est couvrant pour la topologie  f.p.q.c., et pour d\'emontrer
qu'il existe un morphisme $\,h$ comme \c ci-dessus, il suffira donc de v\'erifier que
les deux morphismes 
$h^i_{V_2}=h_V\circ \rm pr_i\,:\,V_2=V\times_X V\lgr Y$ $(i=1,2)$
co\"incident.

Il r\'esulte de la commutativit\'e du premier diagramme de l'\'enonc\'e et 
de l'\'equivariance des morphismes $f$ et $g$ que $g\circ h_V=f\circ v$.
On a donc d\'ej\`a $g\circ h^1_{V_2}=g\circ h^2_{V_2}$. 
Soit $Z^0$ un ouvert dense de $Z$ au-dessus duquel $g$ est une immersion ;
lorsque $T$ est un $Z$-sch\'ema, notons $T^0$ l'ouvert $T\times_Z Z^0$
de $T$. Il r\'esulte alors du fait que $g^0$ est une immersion que  
$ (h^1_{V_2})^0= (h^2_{V_2})^0$. 

Le $X$-sch\'ema $V$ s'identifie (via $(\gamma , u)\mapsto (\gamma, \gamma u)$)
\`a un ouvert de $\Gamma\times_{{\bb Z}_p} X$, et le sch\'ema $V_2$
s'identifie alors \`a un ouvert de 
$\Gamma\times_{{\bb Z}_p}\Gamma\times_{{\bb Z}_p}X$.
En particulier, $V_2$ est r\'eduit et l'ouvert $V^0_2$ est dense dans $V_2$.
De l'identit\'e $ (h^1_{V_2})^0= (h^2_{V_2})^0$ pr\'ec\'edemment obtenue 
il r\'esulte alors que
$ h^1_{V_2}= h^2_{V_2}$ (voir par exemple [Ha, II ex. 4.2]).

Le morphisme $h$ que nous venons d'obtenir est $\Gamma$-\'equivariant par 
construction ; le fait que $h$ est propre lorsque
$f$ l' est  est un corollaire facile du crit\`ere 
valuatif de propret\'e (on rappelle que $g$ est s\'epar\'e). $\square$
\vskip2mm
3.2\pointir Pour d\'emontrer le th\'eor\`eme (2.4.2),
nous utiliserons les deux \'etapes suivantes.

$\underline{\hbox{Etape 1}}$ (cf. 3.3) Nous  d\'efinirons un ouvert
$\cal I$-saturant $U_i$ de  $\widetilde{\cal L}_{\leq i}$
et d\'ecrirons explicitement le sch\'ema $U_i$
(plus pr\'ecis\'ement, nous donnerons en fait des \'equations pour un
$({\bb G}_m ^{i+1})_{{\bb Z}_p}$-torseur $T_i$ au-dessus de $U_i$).
Nous d\'emontrerons ensuite que $U_i$
(ou, ce qui revient au m\^eme, $T_i$) est semi-stable en utilisant le lemme suivant.
\th
LEMME 3.2.1
\enonce Soient $X$ un ${\bb Z}_p$-sch\'ema semi-stable et $Y$ un sous-sch\'ema ferm\'e
de la fibre sp\'eciale $\overline X$ de $X$. Supposons que $Y$ soit lisse sur ${\bb F}_p$ et que
l'intersection de $Y$
et du lieu singulier $\overline X^{\rm sing}$ de $\overline X$ soit
un diviseur \`a croisements normaux (r\'eduit) dans $Y$. L'\'eclat\'e $\widetilde X$
de $X$ le long de $Y$ est alors lui aussi semi-stable.
\endth
{\it D\'emonstration} :  Puisque l'\'enonc\'e que l'on veut d\'emontrer
est local pour la topologie \'etale, on peut supposer 
que $\overline X$ est r\'eunion de $r$ diviseurs lisses $D_i\,$ et travailler
sur un voisinage d'un point ferm\'e $y$ de $Y\cap D_1\cap\cdots\cap D_r$. 
On peut aussi  supposer que  
$Y$ est inclus dans $D_r$. On peut alors  
choisir des \'equations locales $(t_i=0)$ de $D_i$   telles
que $\prod_{1\leq i\leq r}t_i=p$
et des  \'equations locales 
$t_{r+1}=\cdots t_{r+m}=0$ de $Y$
(o\`u $m={\rm codim}_{D_1}\,Y$ et les $t_{r+i}$ sont des sections locales de ${\cal O}_X$). 

Le sch\'ema $Y\cap D_1\cap \cdots\cap D_{r-1}$, qui est 
l'intersection dans $Y$  des $r-1$ branches $Y\cap D_i\  (1\leq i\leq r-1)$ 
du diviseur \`a croisement
normaux $\overline X^{\rm sing}\cap Y$,  est  lisse sur ${\bb F}_p$. Il r\'esulte alors de 
([S], IV. D. 2, Proposition 22) que $t_1, \cdots, t_r,t_{r+1},\cdots,t_{r+m}$ font partie
d'un syst\`eme r\'egulier de param\`etres 
($t_1,\cdots,t_n$) 
de l'anneau local  de $X$ en $y$ (ce dernier  est r\'egulier, puisque $X$ est semi-stable). 
Le morphisme 
$$X\lgr {\rm Spec\,}{\bb Z}_p
[t_1,\cdots,t_n]/(t_1\cdots t_r-p)$$  
d\'efini par ce syst\`eme est alors \'etale au voisinage de $y$. Quitte \`a se localiser
encore sur $X$, on peut donc supposer que
$X={\rm Spec\,}{\bb Z}_p[t_1,\cdots,t_n]/(t_1\cdots t_r-p)$ et que $Y$ est le sous-sch\'ema 
ferm\'e de $X$ d'\'equations $t_r=t_{r+1}=\cdots t_{r+m}=0$. 

Nous allons maintenant d\'ecrire explicitement cet \'eclat\'e.

Consid\'erons le sous-sch\'ema ouvert
$U$ de ${\rm Spec\,}{\bb Z}_p
[\lambda,T_1,\cdots,T_n]/(\lambda T_1\cdots T_r-p)$
compl\'ementaire du sous-sch\'ema ferm\'e d'\'equations $T_r=T_{r+1}=\cdots=T_{r+m}=0$.
On munit $U$ d'une action libre de ${\bb G}_{m,{\bb Z}_p}$ en faisant agir $\Lambda\in 
{\bb G}_{m,{\bb Z}_p}$ par la r\`egle 
$$\lambda\mapsto \Lambda^{-1}\lambda\,,\ T_i\mapsto \Lambda T_i\ (i\in [r,r+m])\,,\ 
T_i\mapsto T_i\ (i\not\in [r,r+m]).$$
On d\'efinit d'autre part un morphisme $U\lgr X$ en envoyant $t_i$ sur $\lambda T_i$ lorsque 
$i\in [r,r+m]$ et sur $T_i$ lorsque $i\not\in [r,r+m]$. Ce morphisme est manifestement
${\bb G}_{m,{\bb Z}_p}$-invariant et induit donc un morphisme 
$U/{\bb G}_{m,{\bb Z}_p}\lgr X$. Il r\'esulte du fait que $t_1,\cdots,t_{r+m}$ est 
une suite r\'eguli\`ere que le morphisme $U/{\bb G}_{m,{\bb Z}_p}\lgr X$ n'est autre que
$\widetilde X\lgr X$  et que le ${\bb G}_{m,{\bb Z}_p}$-torseur
$U\lgr U/{\bb G}_{m,{\bb Z}_p}$ est celui d\'efini par le ${\cal O}_{\widetilde X}$-module  
inversible ${\cal O}_{\widetilde X}(E)$ associ\'e au diviseur exceptionnel $E$.
(cf. [F], Ch IV theorem 2.2). Le ${\bb Z}_p$ sch\'ema $U$ est manifestement semi-stable, et 
donc $\widetilde X=U/{\bb G}_{m,{\bb Z}_p}$ aussi.
$\square$

$\underline{\hbox{Etape 2}}$ (cf. 3.4) Nous construirons (pour $i=3,4,5$) des morphismes 
$$R_i : U_i\lgr  {\rm Gr}_{i-2}\ ({\rm cf.}\  2.3)$$   rendant commutatif le diagramme
$$\kern25mm\diagram{
{U_i\otimes{{\bb Q}_p}} & {\hfl{r_{i}\otimes{{\bb Q}_p}}{}{14mm}} &{{\rm Gr}_{i-2}
\otimes{{\bb Q}_p}}
\cr\noalign{\vskip2mm}
\vfl{}{}{12mm}\cr\noalign{\vskip2mm}
{{\cal L}\otimes{{\bb Q}_p}}&&&\smash{\raise14mm\hbox to 0mm{\kern-24mm 
$\scriptstyle \delta_i 
$\hss}}
\cr
\fleche(38,31)\dir(-1,-1)\long{22}
}$$
 (les produits tensoriels sont sur ${{\bb Z}_p}$ ; $\delta_i$ est le
compos\'e ${\rm pr}_i\circ \delta$, cf. 2.3). 
Les morphismes $R_i$ et le morphisme naturel $U_5\lgr \cal L$ induiront 
alors un morphisme $R_{U_5} : U_5\lgr \prod_{0\leq i \leq 3}{\rm Gr}_i$,
se factorisant par le sous-sch\'ema ferm\'e $\cal M$ puisque $\widetilde{\cal L}_{\leq 5}$
est int\`egre.
En remarquant que la fibre g\'en\'erique du morphisme ${\rm pr}_0\  :\ 
{\cal M}\lgr \cal L$ (cf 2.3) est un isomorphisme,
on applique alors le lemme (3.1.1) au diagramme commutatif \'evident.
On voit ainsi que la construction (3.4) 
d\'emontrera la deuxi\`eme partie du th\'eor\`eme (2.4.2).
\vskip2mm
{\it Remarque} : en utilisant d'une autre mani\`ere le lemme (3.1.1), on voit 
facilement que les morphismes $R_i$ admettent eux aussi un prolongement
\'equivariant sur $\widetilde{\cal L}_{\leq i}$.
\vskip2mm
3.3\pointir Dans les calculs qui suivent, la mani\`ere de noter les coordonn\'es sur l'\'eclat\'e
sera calqu\'ee sur celle employ\'ee  dans la preuve du lemme 3.2.1.
\vskip2mm
3.3.0\pointir Soit $U_0$ l'ouvert ($\cal I$-saturant) 
${\cal L}^0\times_{\cal L}\widetilde{\cal L}_{\leq 0}$ de $\widetilde{\cal L}_{\leq 0}$.
Le sch\'ema 
$U_0$, qui est l'\'eclat\'e de l'espace affine ${\bb A}^6_{{\bb Z}_p}$
(de coordonn\'ees $(a^i_j)_{1\leq i\leq j\leq 3}$) le long de l'origine
de sa fibre sp\'eciale, est alors le quotient
par $({\bb G}_m)_{{\bb Z}_p}$ du ${{\bb Z}_p}$-sch\'ema $T_0$ suivant.

Le sch\'ema $T_0$ est le sous-sch\'ema localement ferm\'e de  
l'espace affine ${\bb A}^8_{{\bb Z}_p}$
(de coordonn\'ees $(\lambda_0, P_0, a^i_j[0])_{1\leq i\leq j\leq 3})$ obtenu 
en intersectant le sous-sch\'ema ferm\'e d'\'equation
$\lambda_0 P_0=p$ avec l'ouvert compl\'ementaire du sous-sch\'ema ferm\'e d'\'equations
$(P_0=a^i_j[0]=0\,,\forall i,j)$. L' action de  $({\bb G}_m)_{{\bb Z}_p}$
sur $T_0$ provient par restriction de celle obtenue en faisant op\'erer
$({\bb G}_m)_{{\bb Z}_p}$
sur ${\bb A}^8_{{\bb Z}_p}$ via l'inverse du caract\`ere tautologique sur le premier
facteur et via le caract\`ere tautologique sur les autres.

On v\'erifie sans difficult\'e (directement ou en utilisant le lemme 3.2.1) que 
$U_0$ est semi-stable.
\vskip2mm
3.3.1\pointir Le sous-sch\'ema ferm\'e ${\cal L}^0_{\{3,5,6\}}$ de ${\cal L}^0$ 
a pour \'equations $a^i_j=0\,,\ \forall (i,j)\not =(3,3)$. Le transform\'e strict 
de sa fibre sp\'eciale est donc le sous-sch\'ema ferm\'e de $U_0$ d\'efini par les 
\'equations $P_0=a^i_j[0]=0\,,\ \forall (i,j)\not =(3,3)$. Comme \c ci-dessus, 
l'\'eclat\'e $U_1$ de $U_0$ le long de ce transform\'e strict est semi-stable et
s'identifie au quotient par 
$({\bb G}_m)^2_{{\bb Z}_p}$ du sch\'ema $T_1$ suivant.

Le sch\'ema $T_1$ est le sous-sch\'ema localement ferm\'e de  
l'espace affine ${\bb A}^9_{{\bb Z}_p}$
(de coordonn\'ees 
$(\lambda_0,\lambda_1, P_1, a^i_j[1])_{1\leq i\leq j\leq 3})$ obtenu 
en intersectant le sous-sch\'ema ferm\'e d'\'equation
$\lambda_0\lambda_1 P_1=p$ avec l'ouvert compl\'ementaire 
des sous-sch\'emas ferm\'es\break d'\'equations respectives
$(P_1=a^i_j[1]=0\,,\forall (i,j)\not=(3,3))$ et $\lambda_1=a^3_3[1]=0$.
 L' action de  $({\bb G}_m)^2_{{\bb Z}_p}$
sur $T_1$ provient de l'action sur  ${\bb A}^9_{{\bb Z}_p}$ d\'efinie par
$$\displaylines{\quad(\Lambda_0,\Lambda_1), (\lambda_0,\lambda_1, 
P_1, a^i_j[1]\,((i,j)\not=(3,3)),a^3_3[1])
\mapsto\hfill\cr
 \hfill(\Lambda_0^{-1}\lambda_0,\Lambda_1^{-1}\lambda_1,
\Lambda_0\Lambda_1 P_1, \Lambda_0\Lambda_1a^i_j[1],\Lambda^0a^3_3[1]).\quad}$$

Soit $U^0_1$ l'ouvert compl\'ementaire du sous-sch\'ema ferm\'e d'\'equation $a^2_3[1]=0$.
A partir du paragraphe (3.3.3) nous travaillerons au-dessus de 
cet ouvert, ce qui nous permettra de 
simplifier l'\'ecriture de certains \'eclatements.
Pour pouvoir  appliquer l'argument (3.2) au-dessus des ouverts 
$U_i=U^0_1\times_{\widetilde{\cal L}_{\leq 1}}\widetilde{\cal L}_{\leq i}$, 
nous aurons alors besoin de savoir que
que ceux-\c ci sont $\cal I$-saturants.
Nous allons donc v\'erifier maintenant que l'ouvert $U^0_1$ est $\cal I$-saturant.

Soient 
$g_1\,,\cdots\,,g_5\in {\cal I}({\bb Z}_p) $ les matrices suivantes.
Pour $1\leq i\leq 4$, $g_i=\pmatrix{\gamma_i&0\cr
0&K{}^t\gamma_i^{-1}K}$, o\`u
$$\gamma_1=\pmatrix{ 
1&0&0\cr
1&1&0\cr
0&0&1\cr}\,,\ 
\gamma_2=\pmatrix{
1&0&0\cr
0&1&0\cr
1&0&1\cr}\,,\ 
\gamma_3=\pmatrix{
1&0&0\cr
0&1&0\cr
0&1&1\cr}\,,\ 
\gamma_4=\pmatrix{
1&0&0\cr
1&1&0\cr
1&0&1\cr}\ ;$$
la matrice $g_5$ s'\'ecrit $\pmatrix{{\rm Id}_3&N\cr
0&{\rm Id}_3\cr}$, o\`u $N=\pmatrix{
0&0&0\cr
p&0&0\cr
0&p&0\cr}$. 
On v\'erifie facilement que l'action de ces cinq \'el\'ements de ${\cal I}({\bb Z}_p) $
sur $\widetilde{\cal L}_{\leq 1}$ pr\'eserve $U'_1$ et que les translat\'es de $U_1$
(qui sont respectivement les ouverts compl\'ementaires des ferm\'es d'\'equation
$a^1_3+a^2_3=0$, $a^1_2+a^2_3=0$, $a^2_2+a^2_3=0$,
$a^1_1+a^1_2+a^1_3+a^2_3=0$ et $P_0+a^2_3=0$)
forment avec $U_1$ un recouvrement de $U'_1$. La surjectivit\'e annonc\'ee
 en r\'esulte imm\'ediatement.
\vskip2mm
{\it Remarque} : le sous-sch\'ema ferm\'e de $T_1$ d'\'equation $a^2_3[1]=1$ s'identifie
bien s\^ur \`a $U^0_1$. Nous pr\'efererons cependant conserver l'homog\'en\'eit\'e
des \'equations en n'utilisant pas cette identification.
\vskip2mm
3.3.2\pointir Le sous-sch\'ema ferm\'e ${\cal L}^0_{\{2,4,6\}}$ de ${\cal L}^0$ 
a pour \'equations $p=a^2_2a^3_3-(a^2_3)^2=a^1_i=0\,,\ \forall i$. 
Le transform\'e strict  $(\overline{{\cal L}^0}_{\{2,4,6\}})\,\widetilde{}\kern-3pt{}_{\leq 1}$
de sa fibre sp\'eciale est donc le sous-sch\'ema ferm\'e de $U_1$ d\'efini par les 
\'equations $P_1=a^2_2[1]a^3_3[1]-\lambda_1(a^2_3[1])^2=a^1_i[1]=0\,,\ \forall i$. 

Pour d\'emontrer que l'\'eclat\'e $U_2$ de $U_1$ le long de $(\overline{{\cal
L}^0}_{\{2,4,6\}})\,\widetilde{}\kern-3pt{}_{\leq 1}$ est semi-stable, nous allons maintenant
v\'erifier les hypoth\`eses du lemme 3.2.1 pour le couple 
$((\overline{{\cal L}^0}_{\{2,4,6\}})\,\widetilde{}\kern-3pt{}_{\leq 1},\  U_1)$. 

On voit
ais\'ement que, localement pour la topologie de Zariski sur 
$(\overline{{\cal L}^0}_{\{2,4,6\}})\,\widetilde{}\kern-3pt{}_{\leq 1}$, l'une au moins des deux
quantit\'es $a^2_2[1]$ et $a^3_3[1]$ est inversible.  La matrice
jacobienne  \vskip 5mm
$$\pmatrix{\smash{\raise 3,5mm\hbox{$\matrix{\partial/\partial\lambda_0\cr\noalign{\vskip3mm}
0}$}}&
\smash{\raise 3,5mm\hbox{$\matrix{\partial/\partial\lambda_1\cr\noalign{\vskip3mm} 0}$}}&
\smash{\raise 3,5mm\hbox{$\matrix{\partial/\partial P_1\cr\noalign{\vskip3mm} 1}$}}&
\smash{\raise 3,5mm\hbox{$\matrix{(\partial/\partial a^1_i[1])_{1\leq i\leq
3}\cr\noalign{\vskip3mm} 0}$}}&
\smash{\raise 3,5mm\hbox{$\matrix{\partial/\partial a^2_2[1]\cr\noalign{\vskip3mm} 0}$}}&
\smash{\raise 3,5mm\hbox{$\matrix{\partial/\partial a^2_3[1]\cr\noalign{\vskip3mm} 0}$}}&
\smash{\raise 3,5mm\hbox{$\matrix{\partial/\partial a^3_3[1]\cr\noalign{\vskip3mm} 0}$}}\cr
0&0&0&0&a^3_3[1]&-2\lambda_1 a^2_3[1]&a^2_2[1]\cr
0&0&0&{\rm Id}_3&0&0&0
} $$
des \'equations de $(\overline{{\cal L}^0}_{\{2,4,6\}})\,\widetilde{}\kern-3pt{}_{\leq
1}$ est donc de rang maximum le long de $(\overline{{\cal
L}^0}_{\{2,4,6\}})\,\widetilde{}\kern-3pt{}_{\leq 1}$. Ceci d\'emontre d\'ej\`a que le sch\'ema
$(\overline{{\cal L}^0}_{\{2,4,6\}})\,\widetilde{}\kern-3pt{}_{\leq 1}$ est lisse sur
${\bb F}_p$. La matrice extraite obtenue en oubliant ses deux premi\`eres colonnes est elle
aussi de rang maximum, ce qui d\'emontre aussi que 
$\overline U_1^{\rm sing}\cap (\overline{{\cal
L}^0}_{\{2,4,6\}})\,\widetilde{}\kern-3pt{}_{\leq 1}$
(qui a pour \'equation $\lambda_0\lambda_1=0$ dans $(\overline{{\cal
L}^0}_{\{2,4,6\}})\,\widetilde{}\kern-3pt{}_{\leq 1}$) est un diviseur \`a coisements normaux
dans $(\overline{{\cal L}^0}_{\{2,4,6\}})\,\widetilde{}\kern-3pt{}_{\leq 1}$.

 Nous allons
maintenant donner une description explicite du sch\'ema $U_2$, comme quotient par 
$({\bb G}_m)^3_{{\bb Z}_p}$ du sch\'ema $T_2$ suivant.

Le sch\'ema $T_2$ est le sous-sch\'ema localement ferm\'e de  
l'espace affine ${\bb A}^{11}_{{\bb Z}_p}$
(de coordonn\'ees 
$(\lambda_0,\lambda_1,\lambda_2, P_2, \delta^1_1[2],
a^i_j[2]_{1\leq i\leq j\leq 3})$ obtenu 
en intersectant le sous-sch\'ema ferm\'e d'\'equations
$$\displaylines{
\kern 20mm\lambda_0\lambda_1\lambda_2 P_2=p\hfill\cr
\kern 20mm\lambda_2\delta^1_1[2]=a^2_2[2]a^3_3[2]-\lambda_1(a^2_3[2])^2\hfill}$$ 
avec l'ouvert compl\'ementaire 
des sous-sch\'emas ferm\'es d'\'equations respectives
$$\displaylines{
\kern 20mm P_2=\delta^1_1[2]=a^1_i=0\,,\ \forall i\hfill\cr  
\kern 20mm\lambda_2=a^2_2[2]=a^2_3[2]=0\hfill\cr 
\kern 20mm\lambda_1=a^3_3[2]=0.\hfill
}$$
 L' action de  $({\bb G}_m)^3_{{\bb Z}_p}$
sur $T_2$ provient de l'action sur  ${\bb A}^{11}_{{\bb Z}_p}$ d\'efinie par
$$\displaylines{\quad(\Lambda_0,\Lambda_1,\Lambda_2)\mapsto\hfill\cr
\hfill{\rm diag}\,
(\Lambda_0^{-1}\!\!,\,\Lambda_1^{-1}\!\!,\,\Lambda_2^{-1}\!\!,\,
\Lambda_0\Lambda_1\Lambda_2\,,\,
(\Lambda_0\Lambda_1)^2\Lambda_2\,,\,(\Lambda_0\Lambda_1\Lambda_2) . {\rm Id}_3
\,,\,(\Lambda_0\Lambda_1).{\rm Id_2}\,,\,
\Lambda_0)\quad}$$
(la notation ${\rm diag}\,(x_1,\cdots,x_n)$ d\'esigne la matrice diagonale
d'\'el\'ements $x_1,\cdots,x_n$ ;
on  range les coordonn\'ees
$a^i_j[2]$ par  $i$ croissants).
\vskip2mm
3.3.3\pointir Le sous-sch\'ema ferm\'e ${\cal L}^0_{\{ 2,3,6\}}$ de ${\cal L}^0$
a pour \'equations $a^1_i=0\,,\ \forall i$.
 Le transform\'e strict de sa fibre 
sp\'eciale a donc pour \'equations $P_2=a^1_i[2]=0\,,\ \forall i$.

Soit $\delta^j_i$ le cofacteur de $a^i_j$ dans la matrice $A$  (cf. 2.1).
Le sous-sch\'ema ferm\'e ${\cal L}^0_{\{ 1,5,6\}}$ de ${\cal L}^0$
a  pour \'equations $\delta^j_i=0\,,\ \forall i,j$. On laisse au lecteur
le soin de v\'erifier que le transform\'e strict de sa fibre sp\'eciale a alors
pour \'equations $P_2=\delta^j_i[2]=0
\,,\,\forall i,j$, o\`u
$$\eqalign{\delta^1_2[2]=-&\left|\matrix{a^1_2[2]&
\lambda_1 a^1_3[2]\cr
a^3_2[2]&\hfill a^3_3[2]\cr}\right|\,,\kern 9pt
\delta^1_3[2]=\left|\matrix{a^2_1[2]&\hphantom{\lambda_2}a^2_2[2]\cr
a^3_1[2]& \hfill a^3_2[2]\cr}\right|
\,,\ \delta^2_2[2]=\left|\matrix{a^1_1[2]&\lambda_1\lambda_2 a^1_3[2]\cr
a^3_1[2]&\hfill a^3_3[2]\cr}\right|\,,\ 
\cr
\delta^2_3[2]=-&\left|\matrix{a^1_1[2]&\lambda_2 a^1_2[2]\cr
a^3_1[2]&\hfill a^3_2[2]\cr}\right|\ {\rm et}\ 
\delta^3_3[2]=\left|\matrix{a^1_1[2]&\lambda_2 a^1_2[2]\cr
a^1_2[2]&\hfill a^2_2[2]\cr}\right|}$$
(avec les conventions $a^i_j[2]=a^j_i[2]$ 
et $\delta^i_j[2]=\delta^j_i[2]$ pour $i\geq j$).

L'intersection dans $U_2$
des transform\'es stricts
de  $\overline{\cal L}^0_{\{ 2,3,6\}}$ et $\overline{\cal L}^0_{\{ 1,5,6\}}$
est vide (cf. la premi\`ere des conditions ouvertes  intervenant dans la d\'efinition 
de $T_2$). Pour d\'emontrer que l'\'eclat\'e de $U_2$ le
long de leur r\'eunion est semi-stable, il suffit donc de v\'erifier que ces deux transform\'es
stricts v\'erifient les hypoth\`eses du lemme 3.2.1. C'est ce ce que nous allons faire maintenant.

Une matrice extraite (en n'oubliant que des colonnes) de la matrice jacobienne des \'equations du
transform\'e strict de 
$\overline{{\cal L}^0}_{\{ 2,3,6\}}$ est 
  \vskip 5mm
$$\pmatrix{
\smash{\raise 3,5mm\hbox{$\matrix{\partial/\partial P_2\cr\noalign{\vskip3mm} 0}$}}&
\smash{\raise 3,5mm\hbox{$\matrix{\partial/\partial \delta^1_1[2]\cr\noalign{\vskip3mm} \lambda_2}$}}&
\smash{\raise 3,5mm\hbox{$\matrix{\partial/\partial a^1_1[2]\cr\noalign{\vskip3mm} 0}$}}&
\smash{\raise 3,5mm\hbox{$\matrix{\partial/\partial a^1_2[2]\cr\noalign{\vskip3mm} 0}$}}&
\smash{\raise 3,5mm\hbox{$\matrix{\partial/\partial a^1_3[2]\cr\noalign{\vskip3mm} 0}$}}&
\smash{\raise 3,5mm\hbox{$\matrix{\partial/\partial a^2_2[2]\cr\noalign{\vskip3mm} -a^3_3[2]}$}}&
\smash{\raise 3,5mm\hbox{$\matrix{\partial/\partial a^3_3[2]\cr\noalign{\vskip3mm} -a^2_2[2]}$}}\cr
1&0&0&0&0&0&0\cr
0&0&1&0&0&0&0\cr
0&0&0&1&0&0&0\cr
0&0&0&0&1&0&0
} $$ 
et, de m\^eme que dans le paragraphe (3.3.2), la proposition \`a v\'erifier r\'esulte du fait
que, localement pour la topologie de Zariski sur $(\overline{{\cal
L}^0}_{\{2,3,6\}})\,\widetilde{}\kern-3pt{}_{\leq 2}$, l'une au moins des trois quantit\'es
$\lambda_2$, $a^2_2[2]$, $a^3_3[2]$ est inversible.

On v\'erifie facilement que
$$
\eqalign{&\lambda_2a^1_3[2]\delta^1_1[2]+a^2_3[2]\delta^1_2[2]
+a^3_3[2]\delta^1_3[2]=0\,,\cr 
&\lambda_2a^1_2[2]\delta^3_1[2]+a^2_2[2]\delta^3_2[2]
+a^3_2[2]\delta^3_3[2]=0\cr
 {\rm et}\quad 
&\lambda_2a^1_3[2]\delta^2_1[2]+a^2_3[2]\delta^2_2[2]
+a^3_3[2]\delta^2_3[2]=0\,.} $$

Au-dessus de l'ouvert $U^0_1$
les quantit\'es $\delta^1_2[2]$, $\delta^3_3[2]$ et $\delta^2_2[2]$
s'expriment alors en fonction de 
$\delta^1_1[2]$, $\delta^1_3[2]$ et $\delta^2_3[2]$. Les \'equations d\'efinissant
le transform\'e
strict de $\overline {\cal L}^0_{\{ 1,5,6\}}$ dans $U^0_2=U_2\times_{U_1}U^0_1$
se r\'esument donc \`a $P_2=\delta^1_1[2]=\delta^1_3[2]=\delta^2_3[2]=0$.
La matrice  \vskip 9mm
$$\pmatrix{
\smash{\raise 3,5mm\hbox{$\matrix{\partial/\partial a^1_1[2]\cr\noalign{\vskip3mm} 0}$}}&
\smash{\raise 3,5mm\hbox{$\matrix{\partial/\partial a^1_2[2]\cr\noalign{\vskip3mm} 0}$}}&
\smash{\raise 3,5mm\hbox{$\matrix{\partial/\partial a^2_2[2]\cr\noalign{\vskip3mm} -a^3_3[2]}$}}&
\smash{\raise 3,5mm\hbox{$\matrix{\partial/\partial \lambda_1\cr\noalign{\vskip3mm}
(a^2_3[2])^2}$}}&
\smash{\raise 3,5mm\hbox{$\matrix{\partial/\partial P_2\cr\noalign{\vskip3mm}
0}$}}&
\smash{\raise 3,5mm\hbox{$\matrix{\partial/\partial \delta^1_1[2]\cr\noalign{\vskip3mm}
\lambda_2}$}}&\cr  
0&a^2_3[2]&-a^1_3[2]&0&0&0\cr
-a^2_3[2]&\lambda_2 a^1_3[2]&0&0&0&0\cr
0&0&0&0&1&0\cr
0&0&0&0&0&1
} $$ 
\vskip2mm
est alors une matrice extraite (en n'oubliant que des colonnes) de la matrice jacobienne des
\'equations de 
$(\overline{{\cal L}^0}_{\{1,4,5\}})\,\widetilde{}\kern-3pt{}_{\leq 2}$ (vu comme quotient par
$({\bb G}_{m,{{\bb Z}_p}})^3$ d'un sous-sch\'ema localement ferm\'e de l'espace affine ambiant
${\bb A}^{11}_{{\bb Z}_p}\supset T_2\supset T^0_2=T_2\times_{U_1}U^0_1$). Il r\'esulte d\'ej\`a du
fait que $a^2_3[2])^2$ est inversible que $(\overline{{\cal
L}^0}_{\{1,4,5\}})\,\widetilde{}\kern-3pt{}_{\leq 2}$ est lisse sur ${\bb F}_p$. En utilisant le
fait que $T^0_2$ est r\'eunion des ouverts compl\'ementaires des ferm\'es d'\'equations
respectives $\lambda_1=0$ et $a^3_3[2]=0$, on voit que l'\'equation
$\lambda_0\lambda_1\lambda_2=0$ d\'efinit un diviseur \`a croisements normaux dans
$(\overline{{\cal L}^0}_{\{1,4,5\}})\,\widetilde{}\kern-3pt{}_{\leq 2}$, ce qui d\'emontre que
$(\overline{{\cal
L}^0}_{\{1,4,5\}})\,\widetilde{}\kern-3pt{}_{\leq 2}$ v\'erifie lui aussi les hypoth\`eses du
lemme 3.2.1.

Nous allons maintenant donner une description explicite de l'\'eclat\'e $U_3$ de
$U^0_2$ le long des transform\'es stricts de 
$\overline{{\cal L}^0}_{\{ 2,3,6\}}$ et $\overline{{\cal L}^0}_{\{ 1,5,6\}}$, comme quotient par
$({\bb G}_m^5)_{{\bb Z}_p}$ du ${\bb Z}_p$-sch\'ema $T_3$ suivant.

Le sch\'ema $T_3$ est le sous-sch\'ema localement ferm\'e de  
l'espace affine ${\bb A}^{15}_{{\bb Z}_p}$
(de coordonn\'ees 
$(\lambda_0,\lambda_1,\lambda_2,\mu_3,\nu_3, P_3, \delta^1_1[3],
\delta^1_3[3],\delta^2_3[3],
a^i_j[3]_{1\leq i\leq j\leq 3})$ obtenu 
en intersectant le sous-sch\'ema ferm\'e  d'\'equations
$$\displaylines{\kern20mm\lambda_0\lambda_1\lambda_2\mu_3\nu_3 P_3=p\hfill\cr
\kern 20mm\lambda_2\nu_3\delta^1_1[3]=a^2_2[3]a^3_3[3]-\lambda_1(a^2_3[3])^2\hfill\cr
\kern 20mm\nu_3\delta^1_3[3]=a^1_2[3]a^2_3[3]-a^1_3[3]a^2_2[3]\hfill\cr
\kern 20mm\nu_3\delta^2_3[3]=-(a^1_1[3]a^2_3[3]-\lambda_2\mu_3a^1_2[3]a^1_3[3])\hfill}$$ 
avec l'ouvert compl\'ementaire 
des sous-sch\'emas ferm\'es d'\'equations respectives
$$\displaylines{\kern 20mm P_3=\delta^1_1[3]=\delta^1_3[3]=\delta^2_3[3]=0\hfill\cr
\kern 20mm\nu_3P_3=a^1_i[3]=0 \,,\  \forall i\hfill\cr
\kern 20mm\mu_3=\nu_3\delta^1_1[3]=0\quad{\rm (cf.\ 3.3.2) }\hfill\cr
\kern 20mm a^2_3[3]=0\hfill\cr 
\kern 20mm\hbox to 0mm{\hss et\quad } 
\lambda_1=a^3_3[3]=0.\hfill}$$
On laisse au lecteur le soin d'expliciter l'action de 
$({\bb G}^5_m)_{{\bb Z}_p}$.
\vskip2mm
3.3.4\pointir Le sous-sch\'ema ferm\'e ${\cal L}^0_{\{ 1,3,5\}}$ de ${\cal L}^0$
a pour \'equations $\delta^i_3= 0\,,\  \forall i$. Le transform\'e strict de
$\overline{{\cal L}^0}_{\{ 1,3,5\}}$ a alors pour \'equations
$P_3=\delta^1_3[3]=\delta^2_3[3]=0$ (on rappelle que 
$\delta^3_3[2]$ s'exprime en fonction de $\delta^1_3[2]$ et $\delta^2_3[2]$,
cf. 3.3.3). De m\^eme que dans le paragraphe 3.3.3, on v\'erifie que ce transform\'e strict
(vu comme sous-sch\'ema ferm\'e de $U_3$) v\'erifie les hypoth\`eses du lemme 3.2.1 (en
n'oubliant que des colonnes, on peut en fait extraire  de la matrice jacobienne de ses
\'equations une matrice d'une forme analogue
\`a celle qu'on a \'ecrite lors de la v\'erification des hypoth\`eses du lemme 3.2.1 pour
$(\overline{{\cal L}^0}_{\{1,4,5\}})\,\widetilde{}\kern-3pt{}_{\leq 2}$).

L'\'eclat\'e $U_4$ de
$U_3$ le long du transform\'e strict de 
 $\overline{{\cal L}^0}_{\{ 1,3,5\}}$
se d\'ecrit explicitement comme le quotient par $({\bb G}_m^6)_{{\bb Z}_p}$ du ${\bb
Z}_p$-sch\'ema
$T_4$ suivant.

Le sch\'ema $T_4$ est le sous-sch\'ema localement ferm\'e de  
l'espace affine ${\bb A}^{16}_{{\bb Z}_p}$
(de coordonn\'ees 
$(\lambda_0,\lambda_1,\lambda_2,\mu_3,\nu_3,\lambda_4, P_4, \delta^1_1[4],
\delta^1_3[4],\delta^2_3[4],
a^i_j[4]_{1\leq i\leq j\leq 3})$ obtenu 
en intersectant le sous-sch\'ema ferm\'e d'\'equations
$$\displaylines{\kern 20mm\lambda_0\lambda_1\lambda_2\mu_3\nu_3\lambda_4 P_4=p\hfill\cr
 \kern 20mm\lambda_2\nu_3\delta^1_1[4]=a^2_2[4]a^3_3[4]-\lambda_1(a^2_3[4])^2\hfill\cr
\kern 20mm\nu_3\lambda_4\delta^1_3[4]=a^1_2[4]a^2_3[4]-a^1_3[4]a^2_2[4]\hfill\cr
\kern 20mm\nu_3\lambda_4\delta^2_3[4]=-(a^1_1[4]a^2_3[4]-\lambda_2\mu_3a^1_2[4]a^1_3[4])\hfill}$$ 
avec l'ouvert compl\'ementaire 
des sous-sch\'emas ferm\'es d'\'equations respectives
$$\displaylines{\kern 20mmP_4=\delta^1_3[4]=\delta^2_3[4]=0\hfill\cr
\kern 20mm\lambda_4=\delta^1_1[4]=0\hfill\cr
\kern 20mm\nu_3\lambda_4P_4=a^1_i[4]=0 \,,\  \forall i\hfill\cr
\kern 20mm\mu_3=\nu_3\delta^1_1[4]=0\quad{\rm (cf.\ 3.3.2) }\hfill\cr
\kern 20mma^2_3[4]=0\hfill\cr 
\kern 20mm\hbox to 0mm{\hss et\quad}\lambda_1=a^3_3[4]=0.\hfill}$$
On laisse au lecteur le soin d'expliciter l'action de 
$({\bb G}^6_m)_{{\bb Z}_p}$.
\vskip2mm
3.3.5\pointir  Le sous-sch\'ema ferm\'e ${\cal L}^0_{\{ 1,2,4\}}$ de ${\cal L}^0$
a pour \'equation d\'et$\, A=0$. Malheureusement,  le 
transform\'e strict de sa fibre sp\'eciale {\it n'est pas} le sous-sch\'ema ferm\'e 
d'\'equation $P_4=\hbox{d\'et}\, A_4=0$, o\`u
$$A_4=\pmatrix{a^1_1[4]&\lambda_2\mu_3a^1_2[4]&\lambda_1\lambda_2\mu_3a^1_3[4]\cr
a^2_1[4]&a^2_2[4]&\lambda_1a^2_3[4]\cr
a^3_1[4]&a^3_2[4]&a^3_3[4]\cr}.$$
En effet, on v\'erifie ais\'ement que 
$-a^2_3[4]\hbox{d\'et}\,A_4=\lambda_2\nu_3^2\lambda_4\Delta$, o\`u
 $$\Delta=\left|\matrix{\delta^1_1[4]&\mu_3\delta^1_2[4]\cr
\delta^3_1[4]&\delta^3_2[4]\cr}\right| {\rm \ et \ }
\delta^1_2[4]=-(a^2_3[4])^{-1}(\lambda_2 a^1_3[4]\delta^1_1[4]
+\lambda_4a^3_3[4]\delta^1_3[4]). $$
Le transform\'e strict de $\overline{{\cal L}^0}_{\{ 1,2,4\}}$ est donc
(\`a priori) un sous-sch\'ema ferm\'e du sous-sch\'ema  ferm\'e 
$(\overline{{\cal L}^0}_{\{1,2,4\}})\,\widetilde{}{\ }^?\kern-10pt{}_{\leq 4}$
de $U_4$ d'\'equations
$P_4=\Delta=0$. Nous allons maintenant d\'emontrer que le couple $((\overline{{\cal
L}^0}_{\{1,2,4\}})\,\widetilde{}{\ }^?\kern-10pt{}_{\leq 4}, U_4)$ satisfait les hypoth\`eses du
lemme 3.2.1. Nous en d\'eduirons ensuite que  $(\overline{{\cal
L}^0}_{\{1,2,4\}})\,\widetilde{}{\ }^?\kern-10pt{}_{\leq 4}$ est en fait le transform\'e strict
$(\overline{{\cal L}^0}_{\{1,2,4\}})\,\widetilde{}\kern -3pt{}_{\leq 4}$ de $\overline{{\cal
L}^0}_{\{1,2,4\}}$ (de sorte que ce transform\'e strict v\'erifiera lui aussi les hypoth\`eses du
lemme 3.2.1).

La matrice
\vskip5mm
$$\pmatrix{
\smash{\raise 3,5mm\hbox {$\matrix{\textstyle\partial\over\textstyle\partial a^1_1[4]\cr\noalign{\vskip3mm} 0}$}}&
\smash{\raise 3,5mm\hbox {$\matrix{\textstyle\partial\over\textstyle\partial a^1_2[4]\cr\noalign{\vskip3mm} 0}$}}&
\smash{\raise 3,5mm\hbox {$\matrix{\textstyle\partial\over\textstyle\partial a^2_2[4]\cr\noalign{\vskip3mm}
-a^3_3[4]}$}}&
\smash{\raise 3,5mm\hbox {$\matrix{\textstyle\partial\over\textstyle\partial a^3_3[4]\cr\noalign{\vskip3mm}
-a^2_2[4]}$}}&
\smash{\raise 3,5mm\hbox {$\matrix{\textstyle\partial\over\textstyle\partial
\lambda_1\cr\noalign{\vskip3mm} (a^2_3[4])^2}$}}&
\smash{\raise 3,5mm\hbox {$\matrix{\textstyle\partial\over\textstyle\partial
\delta^1_1[4]\cr\noalign{\vskip3mm} \lambda_2\nu_3}$}}&
\smash{\raise 3,5mm\hbox {$\matrix{\textstyle\partial\over\textstyle\partial
\delta^2_3[4]\cr\noalign{\vskip3mm}
 0}$}}&
\smash{\raise 3,5mm\hbox {$\matrix{\textstyle\partial\over\textstyle\partial
P_4\cr\noalign{\vskip4,3mm}
 0}$}}\cr \noalign{\vskip1mm}
  0&-a^2_3[4]&a^1_3[4]&0&0&0&0&0\cr
  \kern-3mm a^2_3[4]&\kern-4mm -\lambda_2\mu_3 a^1_3[4]&0&0&0&0&0&0\cr
0&0&0&\kern-3mm \mu_3\lambda_4(\delta^1_3[4])^2&0&\delta^2_3[4]&\delta^1_1[4]&0\cr
0&0&0&0&0&0&0&1
} $$ 
est une matrice extraite (en n'oubliant que des colonnes) de la matrice jacobienne des
\'equations de $(\overline{{\cal L}^0}_{\{1,2,4\}})\,\widetilde{}{\ }^?\kern-10pt{}_{\leq 4}$.
Le long du ferm\'e d'\'equations $\delta^1_1[4]=\delta^2_3[4]=0$, les quantit\'es $\mu_3$,
$\lambda_4$ et $\delta^1_3[4]$ sont inversibles. Par cons\'equent, localement pour la topologie
de Zariski sur $(\overline{{\cal L}^0}_{\{1,2,4\}})\,\widetilde{}{\ }^?\kern-10pt{}_{\leq 4}$,
l'une au moins des entr\'ees de la quatri\`eme ligne de cette matrice est inversible. En
raisonnant sur les trois premi\`eres lignes de cette matrice comme dans les paragraphes 3.3.3 et
3.3.4, ce\c ci d\'emontre que $(\overline{{\cal L}^0}_{\{1,2,4\}})\,\widetilde{}{\
}^?\kern-10pt{}_{\leq 4}$ v\'erifie les hypoth\`eses du lemme 3.2.1.

Le compl\'ementaire dans
$(\overline{{\cal L}^0}_{\{1,2,4\}})\,\widetilde{}{\ }^?\kern-10pt{}_{\leq 4}$ du diviseur
exceptionnel de $U_4$ s'identifie \`a
\`a $\overline{{\cal
L}^0}_{\{1,2,4\}}-\overline{{\cal L}^0}_{\{1,3,5\}}$.  L'intersection de 
$(\overline{{\cal L}^0}_{\{1,2,4\}})\,\widetilde{}{\ }^?\kern-10pt{}_{\leq 4}$ avec le diviseur
exceptionnel de $U_4$ (qui a pour \'equation $\lambda_0\lambda_1\lambda_2\mu_3\nu_3\lambda_4=0$)
est comme on vient de le voir un diviseur \`a croisements normaux dans  $(\overline{{\cal
L}^0}_{\{1,2,4\}})\,\widetilde{}{\ }^?\kern-10pt{}_{\leq 4}$ et est donc de compl\'emntaire
dense. Le sch\'ema $(\overline{{\cal
L}^0}_{\{1,2,4\}})\,\widetilde{}{\ }^?\kern-10pt{}_{\leq 4}$ est donc l'adh\'erence 
sch\'ematique dans $U_4$ de $\overline{{\cal L}^0}_{\{1,2,4\}}-\overline{{\cal
L}^0}_{\{1,3,5\}}$, ce qui d\'emontre que c'est bien le transform\'e strict de $\overline{{\cal
L}^0}_{\{1,2,4\}}$.

On va maintenant d\'ecrire explicitement
l'\'eclat\'e $U_5$ de
$U_4$ le long du transform\'e strict de 
 $\overline{{\cal L}^0}_{\{ 1,2,4\}}$,
comme le quotient par $({\bb G}_m^7)_{{\bb Z}_p}$ du ${\bb Z}_p$-sch\'ema $T_5$
suivant.

Le sch\'ema $T_5$ est le sous-sch\'ema localement ferm\'e de  
l'espace affine ${\bb A}^{18}_{{\bb Z}_p}$
(de coordonn\'ees 
$(\lambda_0,\lambda_1,\lambda_2,\mu_3,\nu_3,\lambda_4,\lambda_5, P_5, 
\Delta[5],\delta^1_1[5],
\delta^1_3[5],\delta^2_3[5],
a^i_j[5]_{1\leq i\leq j\leq 3})$ obtenu 
en intersectant le sous-sch\'ema ferm\'e  d'\'equations
$$\displaylines{\kern 1cm\lambda_0\lambda_1\lambda_2\mu_3\nu_3\lambda_4\lambda_5 P_5=p\hfill\cr
 \kern 1cm\lambda_2\nu_3\delta^1_1[5]=a^2_2[5]a^3_3[5]-\lambda_1(a^2_3[5])^2\hfill\cr
\kern 1cm\nu_3\lambda_4\delta^1_3[5]=a^1_2[5]a^2_3[5]-a^1_3[5]a^2_2[5]\hfill\cr
\kern 1cm\nu_3\lambda_4\delta^2_3[5]=-(a^1_1[5]a^2_3[5]-\lambda_2\mu_3a^1_2[5]a^1_3[5])\hfill\cr
\kern 1cm\lambda_5\Delta[5]=-a^2_3[5]^{-1}(\delta^1_1[5]\delta^3_2[5]+\mu_3\delta^3_1[5]
(a^2_3[5])^{-1}(\lambda_2 a^1_3[5]\delta^1_1[5]+\lambda_4a^3_3[5]\delta^1_3[5]))\hfill}$$ 
avec l'ouvert compl\'ementaire 
des sous-sch\'emas ferm\'es d'\'equations respectives
$$\displaylines{\kern 1cmP_5=\Delta_5=0\hfill\cr
\kern 1cm\lambda_5P_5=\delta^1_3[5]=\delta^2_3[5]=0\hfill\cr
\kern 1cm\lambda_4=\delta^1_1[5]=0\hfill\cr
\kern 1cm\nu_3\lambda_4\lambda_5P_5=a^1_i[5]=0 \,,\  \forall i\hfill\cr
\kern 1cm\mu_3=\nu_3\delta^1_1[5]=0\quad{\rm (cf.\ 3.3.2) }\hfill\cr
\kern 1cm a^2_3[5]=0\hfill\cr 
\kern 1cm\hbox to 0mm{\hss et\quad}
\lambda_1=a^3_3[5]=0.\hfill}$$
On laisse au lecteur le soin d'expliciter l'action de 
$({\bb G}^7_m)_{{\bb Z}_p}$.
\vskip2mm
3.4\pointir Nous allons maintenant construire des morphismes
$R_i\, :\, U_i\lgr {\rm Gr}_{i-2}$ (pour $i=3,4,5$) prolongeant le morphisme
$R^0_i\, :\, U_i\otimes_{{\bb Z}_p}{{\bb Q}_p}\morinj
\widetilde{\cal L}_{\leq i}\otimes_{{\bb Z}_p}{{\bb Q}_p}=\!=
{\cal L}\otimes_{{\bb Z}_p}{{\bb Q}_p}\morinj
{\rm Gr_0}\otimes_{{\bb Z}_p}{{\bb Q}_p}=\!=
{\rm Gr}_{i-2}\otimes_{{\bb Z}_p}{{\bb Q}_p}$ (cf. 2.2.2).
\vskip2mm
3.4.1\pointir Construction du morphisme $R_3$.

On munit le ${\bb Z}_p$-module $V_1$ de 
la base $(\Pi^{-1}e_1,\cdots,\Pi^{-1}e_6)$ (cf. 1.1). Le morphisme $R^0_3$
est donc d\'efini par l'image de la matrice $\Pi\,.\pmatrix{A\cr K}$.
Consid\'erons la matrice
$$M=\pmatrix{\nu_3P_3&0&0\cr
a^1_1[3]&\lambda_0\lambda_1\lambda_2\mu_3a^1_2[3]&
\lambda_0\lambda_1\lambda_2\mu_3a^1_3[3]\cr
a^2_1[3]&\lambda_0\lambda_1a^2_2[3]&\lambda_0\lambda_1a^2_3[3]\cr
a^3_1[3]&\lambda_0\lambda_1a^3_2[3]&\lambda_0\lambda_1a^3_3[3]\cr
0&0&1\cr
0&1&0\cr}.$$
La matrice $\Pi\,.\pmatrix{A\cr K}$ ne diff\`ere de $M$
que par sa premi\`ere colonne, qui est le produit de celle de $M$
par $\lambda_0 \lambda_1\lambda_2\mu_3$. Le morphisme $R^0_i$
est donc encore d\'efini par l'image de $M$.

Les mineurs $\nu_3P_3, a^1_1[3], a^2_1[3], a^3_1[3]$ de la matrice $M$
ne s'annulent pas simultan\'ement sur le sch\'ema $T_3$ et l'image de la 
matrice $M$ est donc un facteur direct de rang 3 du 
$\Gamma\,(T_3,{\cal O}_{T_3})$-module $V_1\otimes_{{\bb Z}_p}
\Gamma\,(T_3,{\cal O}_{T_3})$. Ceci d\'efinit un morphisme
$T_3\lgr {\rm Gr_1}$, et on laisse au lecteur le soin de v\'erifier que
celui-\c ci se factorise par le quotient $T_3\lgr U_3$.
\vskip2mm
{\it Remarque} : la construction pr\'ec\'edente peut d\'ej\`a se faire 
sur l'\'eclat\'e $U_{2,5}$ de $U_2$ le long du transform\'e strict de
$\overline{{\cal L}^0}_{\{ 2,3,6\}}$ 
\vskip2mm
3.4.2\pointir Construction du morphisme $R_4$.

On munit de m\^eme le ${\bb Z}_p$-module $V_2$ de 
la base $(\Pi^{-2}e_1,\cdots,\Pi^{-2}e_6)$, de sorte que le
 morphisme $R^0_4$
est  d\'efini par l'image de la matrice $\Pi^2\,.\pmatrix{A\cr K}$.
Consid\'erons la matrice 
$$B^0=\pmatrix{0&p\cr
p&0\cr
a^1_1&a^1_2
\cr
a^2_1&a^2_2\cr
a^3_1&a^3_2\cr
}$$
form\'ee des deux premi\`eres colonnes et des cinq premi\`eres lignes de la matrice
$\Pi^2\,.\pmatrix{A\cr K}$ et notons $W$ le sous ${\bb Z}_p$-module de $V_2$
engendr\'e par les cinq premiers vecteurs de la base \c ci-dessus. L'image de
la matrice $B^0$ d\'efinit un morphisme 
$b^0\, : \, U_4\otimes_{{\bb Z}_p}{{\bb Q}_p}\morinj
{\cal L}^0\otimes_{{\bb Z}_p}{{\bb Q}_p}\lgr
{\rm Gr}\,(2,W)\otimes_{{\bb Z}_p}{{\bb Q}_p}$ et
le dernier vecteur-colonne
de la matrice $\Pi^2\,.\pmatrix{A\cr K}$ 
d\'efinit un morphisme $c_3\, :\, U_4\lgr{\bb P}(V_2\check{}\,)$ qui se
factorise par
l'ouvert compl\'ementaire de ${\bb P}(W\,\check{}\,)$.
La construction de $R_4$
se ram\`ene donc \`a celle d'un morphisme $b\, :\, U_4\lgr {\rm Gr}\,(2,W)$ 
prolongeant $b^0$ 
(on a alors $R_4=b\oplus c_3$, avec un abus de notations \'evident).
A son tour, la construction du morphisme $b$ se ram\`ene \`a celle d'un
morphisme $\beta\, :\, U_4 \lgr {\bb P}(\Lambda^2 W\,\check{}\,)$
prolongeant le morphisme 
$\beta^0\, :\, U_4 \otimes_{{\bb Z}_p}{{\bb Q}_p}\lgr 
{\bb P}(\Lambda^2 W\,\check{}\,)\otimes_{{\bb Z}_p}{{\bb Q}_p}$
induit par les coordonn\'ees de Pl\"ucker (on rappelle que le ${\bb Z}_p$-sch\'ema
$U_4$ est plat par construction).

Le morphisme $\beta^0$ est d\'efini par les coordonn\'ees homog\`enes 
$$\matrix{\hphantom{M_{1,2}=-\lambda_2\mu_3\nu_3\lambda_4P_4^2,}
&\hphantom{M_{1,3}=-\lambda_2\mu_3P_4a^1_1[4],}\cr
m_{1,2}=-p^2,\hfill 
&m_{1,3}=-pa^1_1,\hfill\cr
m_{1,4}=-pa^2_1,\hfill
&m_{1,5}=-pa^3_1,\hfill\cr
m_{2,3}=pa^2_1,\hfill
&m_{2,4}=pa^2_2,\hfill\cr
m_{2,5}=pa^2_3,\hfill
&m_{3,4}=\delta^3_3,\hfill\cr
m_{3,5}=-\delta^2_3\hfill
{\rm et}\kern-2mm&m_{4,5}=\delta^1_3\,.\hfill}$$
Soient  $$\matrix{
M_{1,2}=-\lambda_2\mu_3\nu_3\lambda_4P_4^2,\hfill 
&M_{1,3}=-\lambda_2\mu_3P_4a^1_1[4],\hfill\cr
M_{1,4}=-\lambda_2\mu_3P_4a^2_1[4],\hfill
&M_{1,5}=-\lambda_2\mu_3P_4a^3_1[4],\hfill\cr
M_{2,3}=\lambda_2\mu_3P_4a^2_1[4],\hfill
&M_{2,4}=P_4a^2_2[4],\hfill\cr
M_{2,5}=P_4a^2_3[4],\hfill
&M_{3,4}=\delta^3_3[4],\hfill\cr
M_{3,5}=-\delta^2_3[4]\hfill
{\rm et}\kern-2mm&M_{4,5}=\delta^1_3[4]\hfill}$$
(o\`u l'on pose $\delta^3_3[4]=
-a^3_2[4]^{-1}(\lambda_2\mu_3a^1_2[4]\delta^3_1[4]+a^2_2[4]\delta^3_2[4])$).
On a alors
$$m_{i,j}=(\lambda_0\lambda_1)^2\lambda_2\mu_3\nu_3\lambda_4M_{i,j}\,,$$
et le morphisme $\beta^0$ est donc encore d\'efini par les coordonn\'ees homog\`enes 
$M_{i,j}$. 

Le triplet $(M_{2,5},M_{3,5},M_{4,5})$ ne s'annule jamais sur
$T_4$ ($T_4$ est recouvert par les ouverts compl\'ementaires des trois ferm\'es
d'\'equations respectives $\delta^1_3[4]=0$, $\delta^1_3[4]=0$ et
$P_4=0$ et $a^2_3[4]$ est inversible sur $T_4$). Le 10-uplet
$(M_{i,j})_{i,j}$ d\'efinit donc un morphisme 
$T_4\lgr {\bb P}(\Lambda^2 W\,\check{}\,)$. On v\'erifie ais\'ement que ce
morphisme est invariant sous l'action naturelle de 
$({\bb G}_m)^6_{{\bb Z}_p}$ sur $T_4$. On obtient alors le morphisme
$\beta$ recherch\'e par passage au quotient. 
\vskip2mm
3.4.3\pointir Construction du morphisme $R_5$.

On munit de m\^eme le ${\bb Z}_p$-module $V_3$ de 
la base $(\Pi^{-3}e_1,\cdots,\Pi^{-3}e_6)$, de sorte que le
 morphisme $R^0_5$
est  d\'efini par l'image de la matrice $\Pi^3\,.\pmatrix{A\cr K}
=\pmatrix{pK\cr A}$.

Le ${\bb Z}_p$-sch\'ema $U_5$ est plat (par construction), et
la construction du morphisme $R_5$ se ram\`ene donc \`a celle d'un
morphisme $\rho_5\, :\, U_5 \lgr {\bb P}(\Lambda^3V_3\,\check{}\,)$
prolongeant le morphisme 
$\rho_5^0\, :\, U_5 \otimes_{{\bb Z}_p}{{\bb Q}_p}\lgr 
{\bb P}(\Lambda^3V_3\,\check{}\,)\otimes_{{\bb Z}_p}{{\bb Q}_p}$
induit par les coordonn\'ees de Pl\"ucker.

Le morphisme $R^0_5$ est d\'efini par les coordonn\'ees homog\`enes
$$\matrix{\hphantom{M_{1,2,5}=
-\lambda_1\lambda_2\mu_3^2\lambda_4\lambda_5P_5^2a^2_1[5]\,,\ \ }
&\hphantom{M_{1,2,4}=
-\lambda_1\lambda_2\mu_3^2\lambda_4\lambda_5P_5^2a^1_1[5]\,,\  }\cr
m_{1,2,3}=-p^3\,,\  \hfill
&m_{1,2,4}=-p^2 a^1_1\,,\  \hfill\cr
m_{1,2,5}=p^2 a^2_1\,,\ \hfill
&m_{1,2,6}=
-p^2 a^3_1\,,\ \hfill\cr
m_{1,3,4}=-m_{1,2,5}\,,\hfill
&m_{1,3,5}=
p^2 a^2_2\,,\ \hfill\cr
m_{1,3,6}=
p^2 a^3_2\,,\ \hfill
&m_{1,4,5}=
p\delta^3_3\,,\ \hfill\cr
m_{1,4,6}=
-p\delta^2_3\,,\ \hfill
&m_{1,5,6}=
p\delta^1_3\,,\hfill\cr
m_{2,3,4}=m_{1,2,6}\,,\ \hfill
&m_{2,3,5}=-m_{1,3,6}\,,\ \hfill\cr
m_{2,3,6}=
-p^2 a^3_3\,,\ \hfill
&m_{2,4,5}=-m_{1,4,6}\,,\ \hfill\cr
m_{2,4,6}=
-p\delta^2_2\,,\hfill
&m_{2,5,6}=
p\delta^1_2\,,\ \hfill\cr
m_{3,4,5}=m_{1,5,6}\,,\ \hfill
&m_{3,4,6}=-m_{2,5,6}\,,\ \hfill\cr
m_{3,5,6}=p\delta^1_1\hfill\rm{et}
&m_{4,5,6}=\hbox{d\'et $A$.}\,\hfill}$$ 

Soient  $$\matrix{M_{1,2,3}=
-\lambda_1\lambda_2\mu_3^2\nu_3\lambda_4^2\lambda_5^2P_5^3\,,\  \hfill
&M_{1,2,4}=
-\lambda_1\lambda_2\mu_3^2\lambda_4\lambda_5P_5^2a^1_1[5]\,,\  \hfill\cr
M_{1,2,5}=
-\lambda_1\lambda_2\mu_3^2\lambda_4\lambda_5P_5^2a^2_1[5]\,,\ \hfill
&M_{1,2,6}=
-\lambda_1\lambda_2\mu_3^2\lambda_4\lambda_5P_5^2a^3_1[5]\,,\ \hfill\cr
M_{1,3,4}=-M_{1,2,5}\,,\hfill
&M_{1,3,5}=
\lambda_1\mu_3\lambda_4\lambda_5P_5^2a^2_2[5]\,,\ \hfill\cr
M_{1,3,6}=
\lambda_1\mu_3\lambda_4\lambda_5P_5^2a^3_2[5]\,,\ \hfill
&M_{1,4,5}=
\lambda_1\mu_3\lambda_4P_5\delta^3_3[5]\,,\ \hfill\cr
M_{1,4,6}=
-\lambda_1\mu_3\lambda_4P_5\delta^2_3[5]\,,\ \hfill
&M_{1,5,6}=
\lambda_1\mu_3\lambda_4P_5\delta^1_3[5]\,,\hfill\cr
M_{2,3,4}=M_{1,2,6}\,,\ \hfill
&M_{2,3,5}=-M_{1,3,6}\,,\ \hfill\cr
M_{2,3,6}=
-\mu_3\lambda_4\lambda_5P_5^2a^3_3[5]\,,\ \hfill
&M_{2,4,5}=-M_{1,4,6}\,,\ \hfill\cr
M_{2,4,6}=
-\mu_3P_5\delta^2_2[5]\,,\hfill
&M_{2,5,6}=
\mu_3P_5\delta^1_2[5]\,,\ \hfill\cr
M_{3,4,5}=M_{1,5,6}\,,\ \hfill
&M_{3,4,6}=-M_{2,5,6}\,,\ \hfill\cr
M_{3,5,6}=P_5\delta^1_1[5]\hfill\rm{et}
&M_{4,5,6}=\Delta[5]\,,\hfill}$$
o\`u l'on pose $$\eqalign{&\delta^1_2[5]=
-a^3_2[5]^{-1}(\lambda_2a^1_3[5]\delta^1_1[5]+\lambda_4a^3_3[5]\delta^1_3[5])\,,\cr 
&\delta^2_2[5]=
-a^3_2[5]^{-1}(\lambda_2\mu_3a^3_1[5]\delta^1_2[5]+\lambda_4a^3_3[5]\delta^3_2[5])
\cr
{\rm et\ }
&\delta^3_3[5]=
-a^3_2[5]^{-1}(\lambda_2\mu_3a^1_2[5]\delta^3_1[5]+a^2_2[5]\delta^3_2[5]).}$$
On a $m_{i,j,k}=\lambda_0^3\lambda_1^2\lambda_2^2\mu_3\nu_3^2\lambda_4\lambda_5
M_{i,j,k}\,$,
et le morphisme $\rho^0_5$ est donc encore d\'efini par les coordonn\'ees homog\`enes 
$M_{i,j,k}\,$. 

  On va maintenant v\'erifier que le 20-uplet
$(M_{i,j,k})_{i,j,k}$ ne s'annule jamais sur $T_5$, ce qui
terminera la construction (il est facile de voir que le morphisme 
$T_5\lgr {\bb P}(\Lambda^3 V_3\check{}\,)$ ainsi d\'efini est
invariant sous l'action naturelle de 
$({\bb G}_m)^7_{{\bb Z}_p}$ sur $T_5$ ; le morphisme $\rho_5$ recherch\'e s'obtient
alors par passage au quotient). 

Consid\'erons le sous-sch\'ema ferm\'e $Z$ de $T_5$ d'\'equations
$M_{i,j,k}=0\,,\ \forall i,j,k$ et supposons qu'il soit non vide.
Le long de $Z$ on a $\Delta[5]=0$ et $P_5$ est donc inversible ; l'annulation de 
$M_{3,5,6}$ entra\^ine alors celle de $\delta^1_1[5]$.  En particulier 
$\mu_3$ et $\lambda_4$ sont donc inversibles le long de $Z$, et il r\'esulte alors 
de l'annulation de $ M_{1,4,6}\,,\ M_{1,5,6}\,,\ M_{2,4,6} \,,\ 
M_{2,5,6}\ {\rm et}\ \delta^1_1[5]$ que $\delta^1_3[5]=\delta^2_3[5]=0$
(on rappelle que $T_5$ est recouvert par les ouverts compl\'ementaires des 
deux ferm\'es d'\'equations respectives $\lambda_1=0$ et $a^3_3[5]=0$).
Le long de $Z$, $\lambda_5P_5$ est donc inversible et il r\'esulte alors de l' 
annulation de $M_{1,3,6}$ et $M_{2,3,6}$ que $\lambda_1=a^3_3[5]=0$, ce qui 
est la contradiction cherch\'ee. 
\vskip 5mm
{\bf 4\pointir Application :  construction d'un mod\`ele semi-stable de ${\cal S}(g,N,p)_{{\bb Q}_p}$
pour $g\leq 3$}
\vskip 2mm
Dans ce paragraphe, $g$ est un entier positif inf\'erieur ou \'egal \`a 3 (le cas 
$g=1$ est cependant trivial).

Consid\'erons le produit fibr\'e
$$\widetilde{\cal T}(g,N,p)={\cal T}(g,N,p)\times_{{\bb M}_g}\widetilde {\cal M}_g\,, $$
o\`u le morphisme ${\cal T}(g,N,p)\lgr{{\bb M}_g}$ est le morphisme $f$ (1.3) et o\`u le 
morphisme $\widetilde {\cal M}_g\lgr{{\bb M}_g}$ est la composition du morphisme
$R\ :\ \widetilde {\cal M}_g=\!=\widetilde {\cal L}_g\lgr{{\cal M}_g}$ de la deuxi\`eme partie
du th\'eor\`eme (2.4.2) avec l'inclusion ferm\'ee ${\cal M}_g\morinj {\bb M}_g$ (2.3).
En utilisant le th\'eor\`eme (2.4.2) et le fait que $f$ est lisse (cf. 1.3), on voit 
imm\'ediatement que le morphisme $\widetilde{\cal T}(g,N,p)\lgr{\cal T}(g,N,p)$
est une r\'esolution semi-stable $\cal I$-\'equivariante (cf. 2.4) de ${\cal T}(g,N,p)$.

L'action du sch\'ema en groupes $\cal I$ sur ${\cal T}(g,N,p)$ est libre, et son action sur 
$\widetilde{\cal T}(g,N,p)$ est donc \`a fortiori libre. Consid\'erons le quotient 
$$\widetilde{\cal S}(g,N,p)=\widetilde{\cal T}(g,N,p)/{\cal I}$$
(qui est  \`a priori un espace alg\'ebrique au sens de Knutson [Kn] ; nous verrons
que pour $g=2$, ce quotient est en fait un sch\'ema quasi-projectif). 
Le th\'eor\`eme suivant r\'esulte alors visiblement de la proposition 1.3.2 et du th\'eor\`eme 2.4.2.
\th
TH\'EOR\`EME 4.1
\enonce
1) Localement pour la topologie \'etale, l'espace alg\'ebrique $\widetilde{\cal S}(g,N,p)$
est ${\bb Z}_p$-isomorphe \`a $\widetilde{\cal M}_g$. En particulier, $\widetilde{\cal S}(g,N,p)$
est semi-stable.

2) Le morphisme $\widetilde{\cal S}(g,N,p)\lgr\widetilde{\cal S}(g,N,p)$
est une r\'esolution semi-stable (cf. 2.4). $\square$
\endth
\vskip2mm
{\it Remarques} : 1) Lorsque $g=2$ (resp. $g=3$), la premi\`ere partie de ce th\'eor\`eme revient
\`a dire que $\widetilde{\cal S}(g,N,p)$ est semi-stable et qu'il existe des points (ferm\'es)
de sa fibre sp\'eciale dont tout voisinage \'etale suffisamment petit est r\'eunion de 
quatre (resp. sept) branches lisses. Puisque $\widetilde{\cal S}(2,N,p)$ 
(resp. $\widetilde{\cal S}(3,N,p)$) est de dimension relative 3 (resp. 6) sur ${\bb Z}_p$, 
le nombre 4 (resp. 7) est le maximum possible pour un espace alg\'ebrique semi-stable.

2) A. J. de Jong a d\'ej\`a construit par une autre m\'ethode 
une r\'esolution semi-stable de ${\cal S}(2,N,p)$,
dont la source est un sch\'ema quasi-projectif 
(cf. [dJ2]). Cette r\'esolution semi-stable
co\"incide avec la n\^otre (cf. la remarque 2 suivant notre th\'eor\`eme 2.4.2), et en particulier
$\widetilde{\cal S}(2,N,p)$ est donc bien un sch\'ema quasi-projectif.

3) On peut obtenir de la m\^eme mani\`ere (cf. la remarque 3 suivant le th\'eor\`eme 2.4.2)
un mod\`ele semi-stable de certaines vari\'et\'es de Shimura 
unitaires  dont le groupe est de type ${\rm GU}(2,2)$ 
en la place archim\'edienne. 
\vskip 20mm
{\bf Bibliographie}
\vskip 5mm
\leftskip=\parindent
\parindent=0mm
\parskip=2mm
\leavevmode\llap{[BC] }J.-F. Boutot et H. Carayol\pointir Uniformisation
$p$-adique des courbes de Shimura, {\sl Ast\'erisque} {\bf 196-197} (1991) pp 45-149

\leavevmode\llap{[CN] }C. N. Chai et P. Norman\pointir Bad reduction of the Siegel moduli
scheme of genus two with $\Gamma_0(p)$-level structure, {\sl Amer. Jour of Math.} {\bf 112}
(1990) pp 1003-1071

\leavevmode\llap{[dJ1] }A. J. de Jong\pointir The moduli spaces of principally
polarized abelian varieties with $\Gamma_0(p)$-level structures,
{\sl J. Alg. Geometry} {\bf 2} (1993) pp 667-688

\leavevmode\llap{[dJ2] }A. J. de Jong\pointir Expos\'e \`a Oberwolfach  
(juillet 1992)

\leavevmode\llap{[F] }W. Fulton\pointir Intersection Theory (Springer-Verlag, 1984)

\leavevmode\llap{[GIT] }J. Fogarty et D. Mumford\pointir Geometric Invariant Theory, 
{\sl second enlarged edition, Ergebnisse der Math. und ihrer Grenzgebiete} {\bf 34} 
(1982)

\leavevmode\llap{[GR] }L. Gruson et M. Raynaud\pointir Crit\`eres de platitude et de
projectivit\'e.
 Techniques de "platification" d'un module, {\sl Inventiones Math.} {\bf 13} (1971) pp. 1-89

\leavevmode\llap{[Ha] }R. Hartshorne\pointir Algebraic Geometry, {\sl Graduate Texts
in Math.} {\bf 52} (1977)

\leavevmode\llap{[KL] }D. Kazhdan et G. Lusztig\pointir Schubert varieties and
Poincar\'e duality,  {\sl Proc. of Symp. in pure Math.} {\bf 36} (1980) pp. 185-203

\leavevmode\llap{[Kn] }D. Knutson\pointir Algebraic Spaces, {\sl Springer Lecture Notes in Math.}
{\bf 203} (1971)

\leavevmode\llap{[LS] }V. Lakshmibai et C. S. Seshadri\pointir Geometry of
$G/P$-II, {\sl Proc. Indian Acad. Sci.} {\bf 87 A} (1978) pp. 1-54

\leavevmode\llap{[Me] }W. Messing\pointir The crystals associated to Barsotti-Tate groups : 
with applications to abelian schemes, {\sl Springer Lecture Notes in Math.} {\bf 264} (1972) 

\leavevmode\llap{[MS] }C. Musili et C. S. Seshadi\pointir Standard Monomial Theory, 
{\sl Springer Lecture Notes in Math.} {\bf 264} (S\'eminaire d'Alg\`ebre Paul Dubreil et 
Marie-Paule Malliavin, 1980) pp. 441-476

\leavevmode\llap{[R] }M. Rapoport\pointir On the bad reduction of Shimura varieties, 
{\sl Automorphic forms, Shimura varieties and $L$-functions vol. II (\'ed. par L. Clozel 
et J. S. Milne),
Persp. in Math. } {\bf 11} (1990) pp. 253-321

\leavevmode\llap{[RZ] }M. Rapoport et Th. Zink\pointir Period Spaces for $p$-divisible
Groups, {\sl Annals of Math. Studies} {\bf 141} (Princeton University Press, 1996)

\leavevmode\llap{[T] }J. Tits\pointir Reductive groups over local fields, 
{\sl Automorphic forms, Representations and $L$-functions vol. I 
(\'ed. par A. Borel et W. Casselman), Proc. of Symp. in pure Math.} {\bf 33 part 1} 
(1979) pp. 29-69
\vskip 5mm
\leavevmode\hfill\vbox{\hbox{Alain Genestier}
\hbox{URA 752 du CNRS}
\hbox{}
\hbox{Universit\'e Paris-Sud}
\hbox{Math\'ematique, B\^at. 425}
\hbox{91405 Orsay (France)}
\hbox{}
\hbox{Adresse \'electronique :}
\hbox{\tt Alain.Genestier@math.u-psud.fr}}
\bye